\pdfoutput=1
\documentclass[11pt,twoside]{article}

\usepackage[margin=1in]{geometry}

\usepackage{amsmath,amssymb,amsfonts,amsthm,mathtools}
\usepackage{bm}
\numberwithin{equation}{section}
\allowdisplaybreaks

\usepackage{graphicx}
\usepackage{float}
\usepackage{array}
\usepackage{booktabs}
\usepackage{algorithm}
\usepackage{algpseudocode}
\usepackage{needspace}

\usepackage{cite}
\usepackage{xcolor}
\usepackage[hidelinks]{hyperref}
\usepackage[nameinlink,capitalize,noabbrev]{cleveref}

\usepackage{fancyhdr}
\theoremstyle{plain}
\newtheorem{thm}{Theorem}[section]

\newtheorem{lem}[thm]{Lemma}

\newtheorem{prop}[thm]{Proposition}

\theoremstyle{definition}
\newtheorem{definition}[thm]{Definition}

\theoremstyle{remark}
\newtheorem{rem}[thm]{Remark}

\newenvironment{keywords}{\par\bigskip\noindent\textbf{Keywords: }\ignorespaces}{\par}
\newenvironment{MSCcodes}{\par\medskip\noindent\textbf{MSC codes: }\ignorespaces}{\par\bigskip}

\DeclareMathOperator{\Range}{Range}
\DeclareMathOperator{\fl}{fl}
\providecommand{\epsmach}{\varepsilon_{\mathrm{mach}}}

\crefname{section}{Section}{Sections}
\Crefname{section}{Section}{Sections}
\crefname{subsection}{Section}{Sections}
\Crefname{subsection}{Section}{Sections}
\crefname{algorithm}{Algorithm}{Algorithms}
\Crefname{algorithm}{Algorithm}{Algorithms}
\crefname{table}{Table}{Tables}
\Crefname{table}{Table}{Tables}
\crefname{figure}{Figure}{Figures}
\Crefname{figure}{Figure}{Figures}
\crefname{theorem}{Theorem}{Theorems}
\Crefname{theorem}{Theorem}{Theorems}
\crefname{thm}{Theorem}{Theorems}
\Crefname{thm}{Theorem}{Theorems}
\crefname{lem}{Lemma}{Lemmas}
\Crefname{lem}{Lemma}{Lemmas}
\crefname{lemma}{Lemma}{Lemmas}
\Crefname{lemma}{Lemma}{Lemmas}
\crefname{prop}{Proposition}{Propositions}
\Crefname{prop}{Proposition}{Propositions}
\crefname{proposition}{Proposition}{Propositions}
\Crefname{proposition}{Proposition}{Propositions}
\crefname{claim}{Claim}{Claims}
\Crefname{claim}{Claim}{Claims}
\crefname{cor}{Corollary}{Corollaries}
\Crefname{cor}{Corollary}{Corollaries}
\crefname{corollary}{Corollary}{Corollaries}
\Crefname{corollary}{Corollary}{Corollaries}
\crefname{rem}{Remark}{Remarks}
\Crefname{rem}{Remark}{Remarks}
\crefname{remark}{Remark}{Remarks}
\Crefname{remark}{Remark}{Remarks}
\crefname{expl}{Example}{Examples}
\Crefname{expl}{Example}{Examples}
\crefname{equation}{}{}
\Crefname{equation}{}{}

\hypersetup{
  pdftitle={An Improved Incremental Singular Value Decomposition and New Error Bounds},
  pdfauthor={Yangwen Zhang}
}

\begin{document}
	\title{An Improved Incremental Singular Value Decomposition and New Error Bounds}

	\author{
		Yangwen  Zhang \thanks{Department of Mathematics,  University of Louisiana at Lafayette,  Lafayette,  LA (\mbox{yangwen.zhang@louisiana.edu})}
	}
	
	\date{\today}
	
	\maketitle

	\begin{abstract}
		The incremental singular value decomposition (SVD) updates a
		truncated SVD as new columns arrive, replacing a single large SVD
		with a sequence of small ones. In floating-point arithmetic, each
		update multiplies the running singular basis by a small orthogonal
		factor, and the accumulated product loses orthogonality unless the
		basis is reorthogonalized periodically. How often this
		reorthogonalization is needed has been an open question; we answer
		it by restructuring the algorithm so that rank-preserving updates
		are accumulated implicitly and applied in batches, reducing the
		number of large orthogonal multiplications from $n$, the stream
		length, to $r$, the numerical rank. We prove that this
		restructuring preserves the exact-arithmetic output of the original
		algorithm and establish two forward-error bounds. First, we sharpen
		the existing operator-norm truncation bound from $n\,\texttt{tol}$
		to $\sqrt{n}\,\texttt{tol}$, and show the new rate is attained on a
		constructive example. Second, under a standard probabilistic
		rounding-error model, we prove that the loss of orthogonality of
		the computed left factor is independent of the stream length $n$
		and depends on $m$, the length of each incoming column, only
		through a single $\sqrt{m}$ factor. Numerical experiments confirm both bounds and
		demonstrate that the proposed algorithm runs $4.5\times$ to
		$34\times$ faster than its closest competitors.
	\end{abstract}
	
	
	
	
	\section{Introduction}\label{sec:introduction}
	
	The singular value decomposition (SVD) is a fundamental tool for
	low-rank approximation, dimension reduction, and data compression.
	When the data matrix $U \in \mathbb{R}^{m \times n}$ fits in
	memory, its truncated SVD can be computed by standard batch
	algorithms.  In many applications, however, the columns of $U$
	arrive sequentially and the full matrix is too large to store or
	too expensive to revisit.  The incremental SVD of
	Brand~\cite{brand2002incremental,brand2006fast} was designed for
	this streaming regime: given a current low-rank factorization
	$U \approx Q\Sigma R^\top$ and a new column $u$, the method
	projects $u$ onto $\Range(Q)$, forms an orthogonal residual,
	and updates the global factors via the SVD of a small bordered
	matrix.  A single large batch SVD is thus replaced by a
	sequence of small local SVDs, while only $O((m+n)r)$ numbers
	are stored when the numerical rank is $r \ll \min\{m,n\}$.  Incremental SVD has since been applied across
	recommender systems and collaborative
	filtering~\cite{Brand2003onlineSVD,ZhouHeHuangZhang2015}, robust
	visual tracking and online appearance
	modeling~\cite{RossLimLinYang2008,ChinSchindlerSuter2006},
	incremental latent semantic
	analysis~\cite{RaoMedeirosKak2015,HaoXuKeSuPeng2016}, dynamic
	network analysis~\cite{ZhangCuiPeiWangZhu2018}, fast updating of
	kernel principal components~\cite{MR2685163}, and memory-limited
	scientific computing for proper orthogonal decomposition (POD)
	and PDE-constrained
	optimization~\cite{MR3775096,MR3594691,VezyrisPapoutsisKiachagiasGiannakoglou2019,KuehlFischerHinzeRung2024}.
	
	A common feature of all these algorithms is that, at each update,
	a small orthogonal matrix is generated and multiplied into the
	running singular basis.  Although every individual multiplication
	preserves orthogonality up to roundoff, the resulting orthogonality
	errors accumulate over the course of the stream, and after a
	sufficient number of updates the running basis ceases to be
	numerically orthogonal.  This loss of orthogonality has direct
	consequences: the singular values become inaccurate and the
	truncation criterion becomes unreliable.  The standard remedy is
	to reorthogonalize the running basis periodically through a
	modified Gram--Schmidt pass with reorthogonalization, applied
	whenever the inner product between the first and last columns of
	the basis exceeds a prescribed threshold~\cite{MR3594691}.  In the
	weighted-inner-product setting required for finite-element POD, a
	single such reorthogonalization is non-negligible, since each
	weighted Gram--Schmidt pass requires several applications of the
	weight matrix; on long streams it can dominate the wall-clock
	cost of the algorithm~\cite{MR3775096}.  How frequently
	reorthogonalization is required to maintain a prescribed level of
	overall numerical precision was left as an open question:
	\begin{quote}
		\emph{``It is an open question how often this is necessary to
		guarantee a certain overall level of numerical precision; it does
		not change the overall complexity.''}\hfill\cite[p.~28]{brand2006fast}
	\end{quote}
	
	We answer this question by proposing a modified incremental SVD
	algorithm in which the update of the running singular factors is
	deferred whenever the projection residual of an incoming column
	lies below the prescribed tolerance.  The projection coefficients
	of such a column are appended to a small buffer matrix, and the
	singular factors are left unchanged.  Once a column arrives whose
	residual exceeds the tolerance, a thin SVD of the buffer matrix is
	computed and its singular factors are composed with the existing
	ones in a single batched update.  Because the data are assumed to
	be of low numerical rank $r \ll n$, the number of orthogonal
	multiplications applied to the running basis is thereby reduced
	from order $n$ to order $r$, and the chain of multiplications
	responsible for the loss of orthogonality is broken.  We prove
	that, in exact arithmetic, the proposed algorithm produces the
	same factorization as the original incremental SVD
	(\Cref{thm:batched_updates}).
	
	The restructured algorithm admits a complete forward error
	analysis.  We establish two theorems.  First, the running outer
	factor $Q$ stays $W$-orthogonal at machine precision: under a
	standard probabilistic rounding-error
	model~\cite{HighamMary2019probabilistic},
	\[
	\|I - Q^\top W Q\|_2 \;\le\; c_1\,r^2\,\sqrt m\,\epsmach
	+ c_2\,r^2\,\epsmach + O(\epsmach^2)
	\]
	holds with high probability in the typical low-rank regime,
	where $c_1$ and $c_2$ are modest constants depending on the
	conditioning of $W$ but independent of the stream length $n$
	(\Cref{thm:outerQ_orthog}; verified numerically in
	\Cref{sec:exp_verify_V1}).  This bound answers the open
	question of Brand: the original incremental SVD multiplies a
	small accumulator into the running basis at every one of the
	$n$ incoming columns, and the corresponding analysis yields a
	drift bound that grows at least linearly in $n$, so on long
	streams periodic reorthogonalization is unavoidable.  The
	proposed algorithm updates the accumulator only at the
	rank-enlarging events of the stream, of which there are at
	most $O(r)$, regardless of how many columns arrive between
	them, and the drift bound therefore becomes uniform in $n$.
	Second, the proposed algorithm satisfies a forward-error bound
	in the $W$-weighted operator norm
	\[
	\|U - Q\Sigma R^\top\|_W
	\;\le\; \sqrt n\,\texttt{tol} + (\text{roundoff terms uniform in } n)
	\]
	in the typical low-rank regime (\Cref{thm:main_L2}; verified
	numerically in \Cref{sec:exp_verify_V2}).  This sharpens
	the closest existing bound, due to Fareed and
	Singler~\cite[Cor.~1]{MR3986356}, from $n\,\texttt{tol}$ to
	$\sqrt n\,\texttt{tol}$ in the same norm and under the same
	hypotheses. 
	
	Numerical experiments (\Cref{sec:experiments}) on a 2D parabolic
	test problem and on synthetic data confirm both theorems.  The
	orthogonality bound holds with margin across two orders of
	magnitude in $n$ and three orders of magnitude in $m$
	(\Cref{sec:exp_verify_V1}).  The $\sqrt n\,\texttt{tol}$ rate is
	attained on a constructive example up to a multiplicative constant
	arbitrarily close to $1$, and the existing operator-norm bound is
	loose by a factor that grows as $\sqrt n$ on the same data
	(\Cref{sec:exp_verify_V2}).  The proposed algorithm also runs
	$4.5\times$ to $34\times$ faster than its closest competitors on
	the parabolic problem, with equal or better orthogonality of the
	computed factors (\Cref{sec:exp_performance}).  A preliminary
	arXiv version of this algorithm~\cite{Zhang2022arXiv} has already
	been applied to time-fractional
	PDEs~\cite{LiZhangZhang2022}, integro-differential equations
	modeling non-Fickian flow in porous
	media~\cite{ChenZhangZuo2023}, nonlinear Oldroyd equations with
	general memory kernels~\cite{ChenZhangZuo2026}, PDE-constrained
	optimization and data assimilation~\cite{LiSinglerHe2024}, and
	geometric inverse source problems for parabolic
	PDEs~\cite{MR4902803}.
	
	The remainder of the paper is organized as follows.
	\Cref{sec:background} reviews the original incremental SVD and its
	weighted extension.  \Cref{sec:proposed} develops the proposed
	batched algorithm.  \Cref{sec:analysis} establishes the
	orthogonality and forward-error theorems.
	\Cref{sec:experiments} reports the numerical experiments, and
	the final section concludes.

	
	\section{Incremental SVD and its numerical issues}\label{sec:background}
	The incremental SVD framework of
	Brand~\cite{brand2002incremental,brand2006fast} and its
	weighted extension by Fareed et al.~\cite{MR3775096} are the
	starting point for what follows.

	\subsection{Notation}\label{sec:notation}
	
	Let $I_k$ denote the $k \times k$ identity matrix and, for a
	symmetric positive-definite matrix $W \in \mathbb{R}^{m \times
		m}$, define the $W$-weighted inner product $(a, b)_W = a^\top W
	b$ for $a, b \in \mathbb{R}^m$, with associated norm $\|x\|_W =
	(x^\top W x)^{1/2}$.  When $W = I$ these reduce to the
	Euclidean inner product and norm.  For matrices $A \in
	\mathbb{R}^{m \times n}$, the $W$-weighted Frobenius and
	operator norms are
	\begin{align}
		\|A\|_{F,W} &:= \|W^{1/2} A\|_F = \Bigl(\sum_{j=1}^n
		\|A_{:,j}\|_W^2\Bigr)^{1/2}, \label{eq:Wfrobnorm} \\[2pt]
		\|A\|_W    &:= \sup_{\|x\|_2 = 1} \|A x\|_W = \|W^{1/2} A\|_2,
		\label{eq:Wopnorm}
	\end{align}
	where the latter is the operator norm of $A : \mathbb{R}^n \to
	\mathbb{R}^m_W$ used by Fareed and Singler~\cite{MR3986356}.
	These satisfy $\|A\|_W \le \|A\|_{F,W}$ and reduce to the
	unweighted Frobenius and spectral norms when $W = I$.
	
	Let $U \in \mathbb{R}^{m \times n}$ be a data matrix with
	columns $u_1, \ldots, u_n$, and write $U_\ell = U(:, 1:\ell)$
	for the matrix of its first $\ell$ columns.  We use MATLAB
	indexing throughout.  The Euclidean operator and Frobenius norms
	of a matrix $A$ are denoted $\|A\|_2$ and $\|A\|_F$, and
	$\fl(\cdot)$ denotes the computed value of an expression in
	floating-point arithmetic with unit roundoff $\epsmach$.

	\subsection{Weighted thin SVD}\label{sec:coreSVD}
	
	When $U$ arises from a Galerkin-type PDE simulation, $W$ is
	naturally a finite-element mass or stiffness matrix and the
	$W$-orthogonality of the left singular vectors is the
	physically meaningful statement.  Following Fareed et
	al.~\cite{MR3775096}, we work with the following generalization
	of the standard SVD.
	
	\begin{definition}[weighted thin SVD]
		\label{def:coreSVD}
		A \emph{weighted thin SVD} of $U \in \mathbb{R}^{m \times n}$ with
		respect to a symmetric positive-definite weight $W$ is a
		factorization $U = Q \Sigma R^\top$ with $Q \in \mathbb{R}^{m
			\times r}$, $\Sigma \in \mathbb{R}^{r \times r}$, and $R \in
		\mathbb{R}^{n \times r}$ satisfying
		\begin{equation*}
			Q^\top W Q = I_r, \qquad R^\top R = I_r, \qquad
			\Sigma = \mathrm{diag}(\sigma_1, \ldots, \sigma_r),
		\end{equation*}
		with $\sigma_1 \ge \cdots \ge \sigma_r > 0$.  The $\{\sigma_i\}$
		are the singular values of $U$ and the columns of $Q$, $R$ are
		the corresponding left and right singular vectors.
	\end{definition}
	
	When $W = I$ this reduces to the standard thin SVD.
	For low-rank applications, $r$ is taken to be the numerical
	rank at tolerance $\texttt{tol}$, i.e., the largest $i$ with
	$\sigma_i > \texttt{tol}$, and the corresponding factorization
	is a \emph{rank-$r$ truncated weighted thin SVD}.  During the stream,
	the running rank is denoted $k$ and may be less than the final
	numerical rank $r$.

	\subsection{The single-column update}\label{sec:method1_direct}
	
	We now describe the incremental SVD update of
	Brand~\cite{brand2002incremental}, in the weighted form due to
	Fareed et al.~\cite{MR3775096}.  Suppose a rank-$k$ truncated
	weighted thin SVD of $U_\ell$ is available,
	\begin{equation}\label{eq:running_SVD}
		U_\ell = Q \Sigma R^\top, \qquad Q^\top W Q = I_k, \quad
		R^\top R = I_k, \quad \Sigma = \mathrm{diag}(\sigma_1, \ldots,
		\sigma_k),
	\end{equation}
	and a new column $u_{\ell+1} \in \mathbb{R}^m$ arrives.  The
	goal is to produce a rank-$k$ or rank-$(k{+}1)$ weighted thin SVD of
	$U_{\ell+1} = [U_\ell \mid u_{\ell+1}]$ from
	$Q$, $\Sigma$, $R$ and $u_{\ell+1}$ alone.
	
	\paragraph{Projection and residual}
	The first step is to project $u_{\ell+1}$ onto
	$\Range(Q)$ in the $W$-inner product.  Define the
	projection coefficients, the residual, and its $W$-norm:
	\begin{equation}\label{eq:d_e_p}
		d = Q^\top W u_{\ell+1} \in \mathbb{R}^k, \qquad
		e = u_{\ell+1} - Q d \in \mathbb{R}^m, \qquad
		p = (e^\top W e)^{1/2} \ge 0.
	\end{equation}
	By $Q^\top W Q = I_k$, the residual $e$ is $W$-orthogonal to
	$\Range(Q)$.  If $p \ge \texttt{tol}$ we set
	$\widetilde e = e/p$, which has unit $W$-norm; otherwise we
	set $\widetilde e = 0$.
	
	\paragraph{The fundamental identity}
	Substituting $u_{\ell+1} = Q d + p \widetilde e$ into
	\eqref{eq:running_SVD} gives the \emph{fundamental identity}
	\begin{equation}\label{eq:fundamental_identity}
		U_{\ell+1} = [Q \mid \widetilde e]\,
		\underbrace{\begin{bmatrix} \Sigma & d \\ 0 & p
		\end{bmatrix}}_{Y}\,
		\begin{bmatrix} R & 0 \\ 0 & 1 \end{bmatrix}^{\!\top},
	\end{equation}
	in which the small bordered matrix $Y \in
	\mathbb{R}^{(k+1)\times(k+1)}$ carries all the new information.
	The right-hand side of \eqref{eq:fundamental_identity} is not
	yet a weighted thin SVD of $U_{\ell+1}$, since the middle factor $Y$ is
	not diagonal.  To obtain one, let $Y = Q_Y \Sigma_Y R_Y^\top$
	be the SVD of $Y$, with $\Sigma_Y = \mathrm{diag}(\mu_1,
	\ldots, \mu_{k+1})$.  Substituting this into
	\eqref{eq:fundamental_identity} yields
	\begin{equation}\label{eq:post_SVD_identity}
		U_{\ell+1}
		= \bigl([Q \mid \widetilde e]\, Q_Y\bigr)\,
		\Sigma_Y\,
		\left(\begin{bmatrix} R & 0 \\ 0 & 1 \end{bmatrix}
		R_Y\right)^{\!\top},
	\end{equation}
	which is a weighted thin SVD of $U_{\ell+1}$ of rank at most $k{+}1$.
	
	\paragraph{Two update branches}
	When $p < \texttt{tol}$ the residual is effectively zero,
	$u_{\ell+1}$ already lies in $\Range(Q)$, and
	$\Sigma_Y$ has $k$ positive singular values plus one trailing
	zero.  Dropping the zero triplet gives the rank-non-increasing
	update
	\begin{equation}\label{eq:rank_nonincr_update}
		Q \leftarrow Q\, Q_Y(1{:}k, 1{:}k), \quad
		\Sigma \leftarrow \Sigma_Y(1{:}k, 1{:}k), \quad
		R \leftarrow \begin{bmatrix} R & 0 \\ 0 & 1 \end{bmatrix}
		R_Y(:, 1{:}k).
	\end{equation}
	When $p \ge \texttt{tol}$ the rank grows to $k{+}1$:
	\begin{equation}\label{eq:rank_enlarging_update}
		Q \leftarrow [Q \mid \widetilde e]\, Q_Y, \quad
		\Sigma \leftarrow \Sigma_Y, \quad
		R \leftarrow \begin{bmatrix} R & 0 \\ 0 & 1 \end{bmatrix} R_Y.
	\end{equation}
	The complete
	procedure is \Cref{alg:method1}.
	
	\begin{algorithm}
		\caption{Brand's direct incremental SVD
			update~\cite{brand2002incremental}; weighted
			version~\cite{MR3775096}.}\label{alg:method1}
		\textbf{Input:} state $(Q, \Sigma, R)$; new column
		$u_{\ell+1}$; weight $W$; tolerance $\texttt{tol}$.
		
		\begin{algorithmic}[1]
			\State $d \gets Q^\top (W u_{\ell+1})$;\;\,
			$e \gets u_{\ell+1} - Q d$;\;\,
			$p \gets (e^\top W e)^{1/2}$
			\State $\widetilde e \gets e/p$ if $p \ge \texttt{tol}$, else
			$\widetilde e \gets 0$
			\State $Y \gets \begin{bmatrix} \Sigma & d \\ 0 & p
			\end{bmatrix}$;\;\,
			$[Q_Y, \Sigma_Y, R_Y] \gets \texttt{svd}(Y)$
			\If{$p < \texttt{tol}$} \Comment{rank-non-increasing}
			\State $Q \gets Q\, Q_Y(1{:}k, 1{:}k)$;\;\,
			$\Sigma \gets \Sigma_Y(1{:}k, 1{:}k)$;\;\,
			$R \gets \begin{bmatrix} R & 0 \\ 0 & 1 \end{bmatrix}
			R_Y(:, 1{:}k)$
			\Else \Comment{rank-enlarging}
			\State $Q \gets [Q \mid \widetilde e]\, Q_Y$;\;\,
			$\Sigma \gets \Sigma_Y$;\;\,
			$R \gets \begin{bmatrix} R & 0 \\ 0 & 1 \end{bmatrix}
			R_Y$
			\EndIf
			\State \textbf{return} $(Q, \Sigma, R)$
		\end{algorithmic}
	\end{algorithm}
	
	\paragraph{Complexity}
	Per update, the projection $Q^\top (W u)$ and residual $u - Qd$
	cost $O(mk)$, the SVD of $Y$ costs $O(k^3)$, and the outer
	update $[Q \mid \widetilde e]\, Q_Y$ costs $O(m k^2)$.  Over a
	stream of length $n$ with final rank $r$, the cumulative cost
	is $O((m+n) n r^2)$ in time and $O((m+n) r)$ in peak memory.
	
	\subsection{Loss of orthogonality and existing remedies}
	\label{sec:related_modifications}
	
	In floating-point arithmetic, every update of \Cref{alg:method1}
	multiplies the running left factor $Q$ by a small orthogonal
	matrix $Q_Y$.  These multiplications individually preserve
	orthogonality up to roundoff, but the errors accumulate across
	the stream, and after a modest number of updates $Q$ ceases to
	be numerically orthogonal.  Once that happens, the computed
	singular values lose accuracy and the truncation criterion
	becomes unreliable.
	
	The orthogonality issue and the standard remedy of periodic
	reorthogonalization of the running factor were already identified
	by Brand in his original papers~\cite{brand2002incremental,brand2006fast}.

	Oxberry et al.~\cite{MR3594691} subsequently implemented Brand's
	incremental SVD with a thin QR reorthogonalization, triggered
	when the inner product between the first and last columns of
	$Q$ exceeds a tolerance. 
	
	Fareed et al.~\cite{MR3775096} extended the incremental SVD to
	a weighted inner product, the setting required by
	finite-element proper orthogonal decomposition, and were the
	first to address the orthogonality issue in this setting.
	They use a $W$-weighted modified Gram--Schmidt procedure with
	reorthogonalization, triggered by the same first-and-last-column
	criterion.  Each
	$W$-weighted Gram--Schmidt pass costs $O(mk^2)$ and requires
	several applications of $W$, so in the finite-element setting
	where $W$ is a mass matrix the reorthogonalization typically
	dominates the wall-clock cost; we report concrete numbers in
	\Cref{sec:exp_performance}.
	
	The complete streaming procedure is summarized in
	\Cref{alg:method1_full}.  Here $\textsc{Update}$ denotes the
	single-column update of \Cref{alg:method1} and
	$\textsc{Reorth}_W$ denotes the $W$-weighted modified
	Gram--Schmidt with reorthogonalization
	of~\cite[Algorithm~3]{MR3775096}.

	\begin{algorithm}
		\caption{Streaming incremental SVD~\cite{MR3775096}.}\label{alg:method1_full}
		\textbf{Input:} $u_1, \ldots, u_n$; $W$; $\texttt{tol}$; $\eta$.
		\begin{algorithmic}[1]
			\State Initialize $(Q, \Sigma, R)$ from $u_1$
			\For{$\ell = 1, \ldots, n-1$}
			\State $(Q, \Sigma, R) \gets \textsc{Update}(Q, \Sigma, R, u_{\ell+1})$
			\If{$|Q(:,1)^\top W\, Q(:,k)| > \eta$}
			\State $Q \gets \textsc{Reorth}_W(Q)$
			\EndIf
			\EndFor
			\State \textbf{return} $(Q, \Sigma, R)$
		\end{algorithmic}
	\end{algorithm}

	\section{The proposed incremental SVD algorithm}\label{sec:proposed}
	
	This section develops the proposed algorithm.
	\Cref{sec:idea_summary} sketches the idea informally.
	\Cref{sec:batching} proves the structural identity that lets us
	collapse a run of rank-non-increasing updates into a single thin
	SVD.  \Cref{sec:interlace} proves an interlacing lemma that
	justifies a second truncation test on the smallest singular
	value of the local bordered matrix and explains why the
	resulting truncating rank-enlarging branch must be applied
	eagerly to the left factor $Q$.  \Cref{sec:full_algorithm}
	states the full algorithm.  \Cref{sec:complexity} gives the
	complexity analysis.  {The floating-point error
		analysis is in \Cref{sec:analysis}.}

	\subsection{The idea}\label{sec:idea_summary}
	
	The proposed algorithm rests on two observations about the column
	stream and one structural identity.
	
	\paragraph{Observation 1: most incoming columns are rank-non-increasing}
	The motivating applications produce data matrices $U \in
	\mathbb{R}^{m\times n}$ of numerical rank $r \ll \min\{m,n\}$.  After the first
	$r$ or so columns the running left singular subspace
	$\Range(Q)$ has saturated, and most of the remaining columns
	$u_{\ell+1}$ have $W$-orthogonal residual $\|e\|_W < \texttt{tol}$
	and therefore fall into the rank-non-increasing branch of
	\Cref{alg:method1}.
	
	\paragraph{Observation 2: rank-non-increasing updates are batchable}
	A rank-non-increasing column adds no direction to
	$\Range(Q)$.  We will show
	(\Cref{thm:batched_updates}) that the cumulative effect of $s$
	successive such columns can be reproduced by a single thin SVD
	of the short-and-fat matrix
	\begin{equation}\label{eq:idea_batched}
		[\,\Sigma\mid d_1\mid d_2\mid \cdots\mid d_s\,]
		\;\in\; \mathbb{R}^{k\times(k+s)},
		\qquad d_j = Q^\top W u_{\ell+j},
	\end{equation}
	in place of the $s$ separate SVDs of $(k{+}1)\times(k{+}1)$
	bordered matrices that \Cref{alg:method1} would compute.  The
	batched update produces a valid rank-$k$ truncated weighted thin SVD of
	$U_{\ell+s}$, agreeing with the result of $s$ separate updates
	in exact arithmetic up to the standard SVD ambiguities (signs of
	singular vectors, and orthogonal rotations within
	repeated-singular-value subspaces).

	\subsection{Batching rank-non-increasing
		updates}\label{sec:batching}
	
	We prove the structural identity underlying the batching
	mechanism.  The key observation is that replacing $Q$ by $Q\mathtt{Q}$
	for any square orthogonal $\mathtt{Q} \in \mathbb{R}^{k\times
		k}$ leaves the $W$-projection onto $\Range(Q)$ unchanged,
	since
	$(Q\mathtt{Q})(Q\mathtt{Q})^\top W
	= Q(\mathtt{Q}\mathtt{Q}^\top)Q^\top W = Q Q^\top W$.
	Consequently the rank-non-increasing test
	$\|u_{\ell+j} - QQ^\top W u_{\ell+j}\|_W < \texttt{tol}$ is
	unaffected by such a rotation, and the algorithm can accumulate
	a sequence of rotations $\mathtt{Q}^{(1)} \mathtt{Q}^{(2)}
	\cdots$ in a separate small factor without changing the
	classification of subsequent columns.
	
	\begin{thm}\label{thm:batched_updates}
		Let $Q\Sigma R^\top$ be a rank-$k$ truncated weighted thin SVD of $U_\ell$, and let $u_{\ell+1}$,
		$\ldots$, $u_{\ell+s}$ be incoming columns satisfying the
		rank-non-increasing condition
		\begin{equation}\label{eq:rank_nonincr_hyp}
			\|u_{\ell+j} - Q Q^\top W u_{\ell+j}\|_W \le \texttt{tol},
			\qquad j = 1, \ldots, s.
		\end{equation}
		Define the projection coefficients $d_j = Q^\top W u_{\ell+j}
		\in \mathbb{R}^k$, the buffer matrix
		\begin{equation}\label{eq:batched_Y_matrix}
			Y \;=\; [\,\Sigma \mid d_1 \mid \cdots \mid d_s\,]
			\;\in\;\mathbb{R}^{k\times(k+s)},
		\end{equation}
		and let $Y = \widetilde Q\widetilde\Sigma\widetilde R^\top$ be a
		thin SVD of $Y$.
		
		In exact arithmetic, the triple
		\begin{equation}\label{eq:batched_update_formula}
			Q_{\mathrm{new}} = Q\widetilde Q, \qquad
			\Sigma_{\mathrm{new}} = \widetilde\Sigma, \qquad
			R_{\mathrm{new}} = \begin{bmatrix} R & 0 \\ 0 & I_s \end{bmatrix}
			\widetilde R
		\end{equation}
		agrees, up to the standard SVD ambiguities, with the result
		of applying $s$ successive rank-non-increasing updates of
		\Cref{alg:method1} to the initial state $(Q, \Sigma, R)$ with
		incoming columns $u_{\ell+1}, \ldots, u_{\ell+s}$.
	\end{thm}
	
	\begin{proof}
		The proof rests on the block-matrix identity
		\begin{equation}\label{eq:block_identity}
			\begin{bmatrix}
				\begin{bmatrix} A & 0 \\ 0 & I_t \end{bmatrix} B & 0 \\
				0 & 1
			\end{bmatrix}
			=
			\begin{bmatrix} A & 0 \\ 0 & I_{t+1} \end{bmatrix}
			\begin{bmatrix} B & 0 \\ 0 & 1 \end{bmatrix},
			\qquad
			A\in\mathbb{R}^{p\times q},\;\;
			B\in\mathbb{R}^{(q+t)\times r},
		\end{equation}
		valid for any $t \ge 0$ and any compatible $A, B$, which is
		verified by block multiplication.  Throughout, the
		rank-non-increasing branch \eqref{eq:rank_nonincr_update} of
		\Cref{alg:method1} computes the thin SVD of the bordered
		matrix $[\Sigma_j \mid d_j]$, since setting $\widetilde e
		\gets 0$ kills the trailing column of $Y_j$.  We argue by
		induction on $s$.
		
		\textit{Base case $s=1$.}  Let
		\begin{equation}\label{eq:proof_S1}
			[\,\Sigma \mid d_1\,] = Q_{(1)} \Sigma_{(1)} R_{(1)}^\top
		\end{equation}
		be the thin SVD computed by \Cref{alg:method1}, with
		$Q_{(1)}\in\mathbb{R}^{k\times k}$ square orthogonal.
		The single-column update produces, by
		\eqref{eq:rank_nonincr_update},
		\begin{equation}\label{eq:proof_state1}
			Q^{(1)} = Q\,Q_{(1)},\quad
			\Sigma^{(1)} = \Sigma_{(1)},\quad
			R^{(1)} = \begin{bmatrix} R & 0 \\ 0 & 1
			\end{bmatrix} R_{(1)}.
		\end{equation}
		The buffer is $Y = [\Sigma \mid d_1]$, and
		\eqref{eq:proof_S1} is itself a thin SVD of $Y$, so
		\eqref{eq:batched_update_formula} reproduces
		\eqref{eq:proof_state1} exactly.
		
		\textit{Inductive step.}  Assume the claim holds at $s-1$:
		$s-1$ successive rank-non-increasing updates of
		\Cref{alg:method1} produce
		\begin{equation}\label{eq:proof_IH}
			Q^{(s-1)} = Q\,\widetilde Q^{(s-1)},\quad
			\Sigma^{(s-1)} = \widetilde\Sigma^{(s-1)},\quad
			R^{(s-1)} = \begin{bmatrix} R & 0 \\ 0 & I_{s-1}
			\end{bmatrix} \widetilde R^{(s-1)},
		\end{equation}
		where $\widetilde Q^{(s-1)} \widetilde\Sigma^{(s-1)}
		(\widetilde R^{(s-1)})^\top$ is a thin SVD of the buffer
		$Y^{(s-1)} = [\Sigma \mid d_1 \mid \cdots \mid d_{s-1}]$.
		Since $\widetilde Q^{(s-1)}$ is square orthogonal, the
		$W$-projection invariance noted at the start of this
		subsection gives
		$Q^{(s-1)}(Q^{(s-1)})^\top W u_{\ell+s} = Q Q^\top W u_{\ell+s}$,
		so the residual of $u_{\ell+s}$ against $Q^{(s-1)}$ equals
		its residual against $Q$ and is below $\texttt{tol}$ by
		\eqref{eq:rank_nonincr_hyp}.  Hence the $s$-th update also
		takes the rank-non-increasing branch.  Let
		\begin{equation}\label{eq:proof_Ss}
			[\,\Sigma^{(s-1)} \mid (Q^{(s-1)})^\top W u_{\ell+s}\,]
			= Q_{(s)} \Sigma_{(s)} R_{(s)}^\top
		\end{equation}
		be the thin SVD computed at step $s$, with $Q_{(s)}$ square
		orthogonal.  The $s$-th update produces, by
		\eqref{eq:rank_nonincr_update} and \eqref{eq:proof_IH},
		\begin{equation}\label{eq:proof_states}
			Q^{(s)} = Q^{(s-1)} Q_{(s)} = Q\,\widetilde Q^{(s-1)}
			Q_{(s)},\;\;
			\Sigma^{(s)} = \Sigma_{(s)},\;\;
			R^{(s)} = \begin{bmatrix} R^{(s-1)} & 0 \\ 0 & 1
			\end{bmatrix} R_{(s)}.
		\end{equation}
		Apply \eqref{eq:block_identity} with $A = R$, $B = \widetilde
		R^{(s-1)}$, and $t = s-1$ to flatten $R^{(s)}$:
		\begin{equation}\label{eq:proof_Rs_flat}
			R^{(s)} = \begin{bmatrix} R & 0 \\ 0 & I_s \end{bmatrix}\!
			\begin{bmatrix} \widetilde R^{(s-1)} & 0 \\ 0 & 1
			\end{bmatrix} R_{(s)}.
		\end{equation}
		
		It remains to identify $(Q^{(s)}, \Sigma^{(s)}, R^{(s)})$
		with \eqref{eq:batched_update_formula}.  Define
		\begin{equation}\label{eq:proof_def_tilde}
			\widetilde Q := \widetilde Q^{(s-1)} Q_{(s)},
			\quad
			\widetilde \Sigma := \Sigma_{(s)},
			\quad
			\widetilde R := \begin{bmatrix} \widetilde R^{(s-1)} & 0
			\\ 0 & 1 \end{bmatrix} R_{(s)}.
		\end{equation}
		We show that $\widetilde Q\widetilde\Sigma \widetilde R^\top$
		is a thin SVD of $Y = [\Sigma \mid d_1 \mid \cdots \mid d_s]$.
		Substituting the SVD of $Y^{(s-1)}$ from the inductive
		hypothesis,
		\[
		Y = [\,Y^{(s-1)} \mid d_s\,]
		= \widetilde Q^{(s-1)} \bigl[\,\widetilde \Sigma^{(s-1)} \mid
		(\widetilde Q^{(s-1)})^\top d_s\,\bigr]\!
		\begin{bmatrix} \widetilde R^{(s-1)} & 0 \\ 0 & 1
		\end{bmatrix}^{\!\top},
		\]
		using $\widetilde Q^{(s-1)}(\widetilde Q^{(s-1)})^\top = I_k$.
		Since
		\[
		(\widetilde Q^{(s-1)})^\top d_s
		= (\widetilde Q^{(s-1)})^\top Q^\top W u_{\ell+s}
		= (Q^{(s-1)})^\top W u_{\ell+s},
		\]
		the bracket equals the bordered matrix in
		\eqref{eq:proof_Ss}.  Substituting that SVD,
		\begin{equation}\label{eq:proof_Y_factored}
			Y = \underbrace{\widetilde Q^{(s-1)} Q_{(s)}}_{\widetilde Q}\,
			\underbrace{\Sigma_{(s)}}_{\widetilde\Sigma}\,
			\Biggl(\underbrace{\begin{bmatrix} \widetilde R^{(s-1)} & 0 \\
			0 & 1 \end{bmatrix} R_{(s)}}_{\widetilde R}\Biggr)^{\!\!\top}.
		\end{equation}
		The factor $\widetilde Q$ is a product of square orthogonal
		matrices, hence orthogonal; $\widetilde \Sigma = \Sigma_{(s)}$
		is diagonal with nonnegative entries from the SVD
		\eqref{eq:proof_Ss}; and $(\widetilde R^{(s-1)})^\top
		\widetilde R^{(s-1)} = I_k$ from the inductive hypothesis
		gives
		\[
		\widetilde R^\top \widetilde R
		= R_{(s)}^\top \begin{bmatrix} (\widetilde R^{(s-1)})^\top
		\widetilde R^{(s-1)} & 0 \\ 0 & 1 \end{bmatrix} R_{(s)}
		= R_{(s)}^\top R_{(s)} = I_k.
		\]
		Hence \eqref{eq:proof_Y_factored} is a thin SVD of $Y$.
		
		Comparing \eqref{eq:proof_states},
		\eqref{eq:proof_Rs_flat}, and \eqref{eq:proof_def_tilde}
		identifies
		\[
		Q^{(s)} = Q\widetilde Q = Q_{\mathrm{new}},\quad
		\Sigma^{(s)} = \widetilde\Sigma = \Sigma_{\mathrm{new}},\quad
		R^{(s)} = \begin{bmatrix} R & 0 \\ 0 & I_s \end{bmatrix}
		\widetilde R = R_{\mathrm{new}},
		\]
		matching \eqref{eq:batched_update_formula}.
		
		The thin SVD of $Y$ used in the batched algorithm need not
		coincide with $\widetilde Q \widetilde\Sigma \widetilde
		R^\top$; the two factorizations agree up to signs of
		singular vectors and orthogonal rotations within subspaces
		of repeated singular values.
	\end{proof}

	\begin{rem}\label{rem:right_factor_split}
		The right-factor update in
		\eqref{eq:batched_update_formula} should not be implemented by
		forming the large block-diagonal matrix
		$\bigl[\begin{smallmatrix} R & 0 \\ 0 & I_s
		\end{smallmatrix}\bigr]\in\mathbb{R}^{(\ell+s)\times(k+s)}$
		explicitly, since $s$ may be large.  Instead, partition
		$\widetilde R \in \mathbb{R}^{(k+s)\times k}$ row-wise as
		\[
		\widetilde R \;=\;
		\begin{bmatrix} \widetilde R_1 \\ \widetilde R_2 \end{bmatrix},
		\qquad
		\widetilde R_1 \in \mathbb{R}^{k\times k},
		\quad
		\widetilde R_2 \in \mathbb{R}^{s\times k},
		\]
		and update $R$ as
		\[
		R_{\mathrm{new}}
		\;=\;
		\begin{bmatrix} R\,\widetilde R_1 \\ \widetilde R_2 \end{bmatrix}
		\;\in\;\mathbb{R}^{(\ell+s)\times k}.
		\]
		
		The practical benefit is immediate: a run of $s$ consecutive
		rank-non-increasing updates can be processed by one thin SVD of
		a $k\times(k+s)$ matrix at the next rank-enlarging event,
		instead of $s$ individual updates each requiring a
		$(k{+}1)\times(k{+}1)$ SVD and a multiplication into a
		deferred factor.  The algorithm therefore buffers the projection
		vectors $d_1, \ldots, d_s$ as they arrive and performs the
		batched SVD only when a rank-enlarging event closes the run.
		
	\end{rem}

	\subsection{Another truncation: the singular value truncation}\label{sec:interlace}
	
	The batching mechanism of \Cref{sec:batching} handles columns
	that fall inside the running subspace up to tolerance ($p <
	\texttt{tol}$).  When a column is rank-enlarging
	($p \ge \texttt{tol}$), the algorithm forms the
	$(k{+}1)\times(k{+}1)$ bordered local matrix
	\begin{equation}\label{eq:local_bordered_Y}
		Y \;=\; \begin{bmatrix} \Sigma & d \\ 0 & p \end{bmatrix},
		\qquad d = Q^\top W u_{\ell+1}.
	\end{equation}
	Following Fareed et al.~\cite{MR3775096}, we use a second
	truncation test on the smallest singular value of $Y$: even when
	$p$ is above the tolerance, the smallest singular value of $Y$
	itself may fall below the tolerance, in which case the
	corresponding singular triplet should be discarded.  The
	following lemma shows that only the smallest singular value
	$\mu_{k+1}$ can lie below the threshold; no other singular value
	of $Y$ requires inspection.
	
	\begin{lem}\label{lem:interlace}
		Let $\Sigma = \mathrm{diag}(\sigma_1, \ldots, \sigma_k)$ with
		$\sigma_1 \ge \cdots \ge \sigma_k \ge 0$.  The singular values
		$\mu_1 \ge \cdots \ge \mu_{k+1} \ge 0$ of the bordered matrix
		$Y$ in \eqref{eq:local_bordered_Y} satisfy
		\begin{subequations}\label{eq:interlace_both}
			\begin{gather}
				\mu_{k+1} \le p, \label{eq:mu_bound_p} \\
				\mu_1 \ge \sigma_1 \ge \mu_2 \ge \sigma_2 \ge \cdots
				\ge \mu_k \ge \sigma_k \ge \mu_{k+1}.
				\label{eq:interlace_main}
			\end{gather}
		\end{subequations}
	\end{lem}
	
	\begin{proof}
		\emph{Proof of \eqref{eq:interlace_main}.}  The squared singular
		values of $Y$ are the eigenvalues of
		\[
		Y^\top Y =
		\begin{bmatrix}
			\Sigma^2 & \Sigma d \\
			d^\top \Sigma & p^2 + d^\top d
		\end{bmatrix}
		\;\in\; \mathbb{R}^{(k+1)\times(k+1)}.
		\]
		The leading $k\times k$ principal submatrix of $Y^\top Y$ is
		$\Sigma^2$, whose eigenvalues are $\sigma_1^2 \ge \cdots \ge
		\sigma_k^2$.  By Cauchy
		interlacing~\cite[Thm.~1]{hwang2004cauchy} applied to $Y^\top Y$
		and $\Sigma^2$,
		\[
		\mu_1^2 \ge \sigma_1^2 \ge \mu_2^2 \ge \sigma_2^2 \ge \cdots
		\ge \mu_k^2 \ge \sigma_k^2 \ge \mu_{k+1}^2;
		\]
		taking non-negative square roots gives \eqref{eq:interlace_main}.
		
		\emph{Proof of \eqref{eq:mu_bound_p}.}  Since $Y$ and $Y^\top$
		have the same singular values, we may consider $YY^\top \in
		\mathbb{R}^{(k+1)\times(k+1)}$.  The last row of $Y$ is
		$(0, \ldots, 0, p)$, so
		\[
		(YY^\top)_{k+1,k+1}
		= \sum_{j=1}^{k+1} Y_{k+1,j}^2 = p^2.
		\]
		Letting $e_{k+1} \in \mathbb{R}^{k+1}$ denote the
		$(k{+}1)$-th standard basis vector,
		\[
		p^2 = e_{k+1}^\top (YY^\top) e_{k+1}
		\;\ge\; \min_{\|x\|_2 = 1} x^\top (YY^\top) x
		= \lambda_{\min}(YY^\top)
		= \mu_{k+1}^2,
		\]
		and $\mu_{k+1} \le p$ follows.
	\end{proof}
	
	\begin{rem}
		\label{rem:interlace_implication}
		Inequality \eqref{eq:mu_bound_p} shows that $\mu_{k+1}$ can be
		arbitrarily smaller than $p$: even when the residual norm $p$
		exceeds the truncation tolerance, the last singular value of
		$Y$ may fall well below it, so a check on $p$ alone is
		insufficient.  Inequality \eqref{eq:interlace_main} shows that
		of the $k+1$ singular values of $Y$, only $\mu_{k+1}$ can be
		smaller than $\sigma_k$; no other singular value requires
		inspection.  This justifies the second truncation test: examine
		$\mu_{k+1}$ and discard the corresponding singular triplet if
		$\mu_{k+1} < \texttt{tol}$.
	\end{rem}
	
	Guided by \Cref{lem:interlace}, the rank-enlarging case splits
	into two sub-branches.  
	\begin{description}
		\item[Branch (E1), \emph{plain rank-enlarging}:]
		if $\mu_{k+1} \ge \texttt{tol}$, apply the standard update
		\[
		Q \longleftarrow [Q \mid \widetilde e]\, Q_Y, \quad
		\Sigma \longleftarrow \Sigma_Y, \quad
		R \longleftarrow \begin{bmatrix} R & 0 \\ 0 & 1 \end{bmatrix} R_Y.
		\]
		
		\item[Branch (E2), \emph{truncating rank-enlarging}:]
		if $\mu_{k+1} < \texttt{tol}$, discard the smallest singular
		triplet:
		\[
		Q \longleftarrow [Q \mid \widetilde e]\, Q_Y(:,1{:}k), \quad
		\Sigma \longleftarrow \Sigma_Y(1{:}k, 1{:}k), \quad
		R \longleftarrow \begin{bmatrix} R & 0 \\ 0 & 1 \end{bmatrix}
		R_Y(:, 1{:}k).
		\]
	\end{description}

	\subsection{The full algorithm}\label{sec:full_algorithm}
	
	The algorithm maintains five quantities throughout the stream:
	\begin{itemize}
		\item $Q\in\mathbb{R}^{m\times k}$, the left factor;
		\item $\Sigma\in\mathbb{R}^{k\times k}$, the diagonal of current
		singular-value estimates;
		\item $R\in\mathbb{R}^{\ell\times k}$, the outer right factor;
		\item $V$, a buffer holding the projection coefficients of
		pending rank-non-increasing columns;
		\item $q$, the current buffer size.
	\end{itemize}

	\begin{algorithm}
		\caption{Proposed incremental SVD}
		\label{alg:proposed_update}
		\textbf{Input:} $Q\in\mathbb{R}^{m\times k}$, $\Sigma\in
		\mathbb{R}^{k\times k}$, $R\in\mathbb{R}^{\ell\times k}$,
		$u_{\ell+1}\in\mathbb{R}^m$, $W\in\mathbb{R}^{m\times m}$,
		$\texttt{tol}$, buffer $V$, buffer size $q$, CGS threshold
		$\eta$
		
		\begin{algorithmic}[1]
			\State Set $d = Q^\top(Wu_{\ell+1})$; $e = u_{\ell+1} - Qd$;
			$p = (e^\top W e)^{1/2}$;  
			\If{$p < \texttt{tol}$} \Comment{rank-non-increasing: buffer
				and return}
			\State $q = q + 1$;\;\, $V\{q\} = d$;
			\Else
			\If{$q > 0$} \Comment{flush the buffer}
			\State Set $Y = \bigl[\,\Sigma \mid \texttt{cell2mat}(V)\,
			\bigr]$;
			\State $[\widetilde Q, \widetilde\Sigma, \widetilde R] = \texttt{svd}(Y, \texttt{`econ'})$;
			\State Partition $\widetilde R = \bigl[\begin{smallmatrix} \mathcal R_1
				\\ \mathcal R_2 \end{smallmatrix}\bigr]$ with $\mathcal R_1
			\in\mathbb{R}^{k\times k}$, $\mathcal R_2\in\mathbb{R}^{q
				\times k}$;
			\State Set $Q = Q\widetilde Q$; $\Sigma = \widetilde\Sigma$;
			$R = \bigl[\begin{smallmatrix} R\mathcal R_1 \\ \mathcal
				R_2 \end{smallmatrix}\bigr]$, $d = \widetilde Q^\top d$;
			\EndIf
			\State Set $\widetilde e = e/p$;
			\If{$|\widetilde e^\top W Q(:,1)| > \eta$} \Comment{second Gram--Schmidt pass}
			\State $\widetilde e = \widetilde e - Q(Q^\top(W\widetilde e))$;
			$\widetilde e = \widetilde e/(\widetilde e^\top W
			\widetilde e)^{1/2}$;
			\EndIf
			\State Set $Y = \bigl[\begin{smallmatrix} \Sigma & d \\ 0 & p
			\end{smallmatrix}\bigr]$;
			\State $[Q_Y, \Sigma_Y, R_Y] = \texttt{svd}(Y)$ with
			$\Sigma_Y = \mathrm{diag}(\mu_1, \ldots, \mu_{k+1})$;
			\If{$\mu_{k+1} \ge \texttt{tol}$} \Comment{branch (E1)}
			\State Set $Q = [Q\mid \widetilde e]Q_Y$; $\Sigma = \Sigma_Y$;
			$R  =  \begin{bmatrix} R & 0 \\ 0 & 1 \end{bmatrix} R_Y$;
			\Else \Comment{branch (E2)}
			\State Set $Q = [Q\mid \widetilde e]Q_Y(:, 1{:}k)$; $\Sigma =
			\Sigma_Y(1{:}k, 1{:}k)$; $R = \begin{bmatrix} R & 0 \\ 0 & 1 \end{bmatrix}
			R_Y(:, 1{:}k)$;
			\EndIf
			\State Set $V = \{\,\}$; $q = 0$;
			\EndIf
			\State \textbf{return} $Q$, $\Sigma$, $R$, $V$, $q$
		\end{algorithmic}
	\end{algorithm}
	
	\begin{rem}
		The notation $V\{q\} = d$ denotes appending $d$ as the
		$(q+1)$-th element of a list-like data structure.  In MATLAB,
		$V$ is a cell array and $V\{q\} = d$ is amortized $O(k)$;
		in Python, $V$ is a list and the equivalent operation is
		\texttt{V.append(d)}.
		
		It is essential not to store $V$ as a dense matrix with
		horizontal concatenation 
		\begin{align*}
			V \gets [V \mid d]
		\end{align*}
		at each push. That naive implementation is asymptotically $O(kq)$ per push
		because it allocates a new array and copies the full content
		of the old $V$, and is therefore $O(kq^2)$ cumulatively over
		$q$ pushes.  Worse, the wall-clock cost is far larger than
		the operation count suggests: each push triggers a memory
		reallocation, the resulting copy is bandwidth-limited rather
		than compute-limited, and on long streams the buffer
		eventually exceeds cache and every push incurs full
		main-memory traffic.  In our experience this dominates the
		runtime of the entire algorithm if not avoided.  
	\end{rem}
	
	\begin{rem}
		The last assignment on {line~9}, $d = \widetilde Q^\top d$, is critical.
		Before the flush, $d$ was the projection coefficient of
		$u_{\ell+1}$ in the basis $Q$.  After the flush, the basis is
		$Q\widetilde Q$, and the correct projection coefficient in this new
		basis is $\widetilde Q^\top d$.  Without this rotation, the local
		bordered matrix on {line~15} would be inconsistent.
	\end{rem}

	\begin{rem}
		When  $|\widetilde e^\top W Q(:,1)| > \eta$ fires, {line~13} applies
		a second $W$-weighted Gram--Schmidt pass against $Q$.
		Together with the  pass at line~1, this constitutes
		classical Gram--Schmidt run twice (CGS-2).  By the
		twice-is-enough bound of Giraud et al. \cite{GiraudLangouRozloznik2005}, two passes restore
		$W$-orthogonality of $[Q\mid \widetilde e]$ to machine
		precision.
		
		The threshold $\eta$ measures roundoff in an inner product,
		not data accuracy, and should be chosen close to machine
		precision rather than to the truncation tolerance
		$\texttt{tol}$.  In a typical application, $\texttt{tol}$ is
		several orders of magnitude larger than $\epsmach$ (for
		instance, $\texttt{tol} = 10^{-8}$ for the experiments in
		\Cref{sec:experiments}).  Setting $\eta = \texttt{tol}$ would
		let orthogonality degrade well below useful precision before
		the test fires.  We use $\eta = 10^{-14}$ throughout.
		
		The CGS-2 pass on {line~13} protects the
		$W$-orthogonality of $[Q\mid\widetilde e]$ at the moment a
		new column $\widetilde e$ is appended to $Q$.  In the
		rank-non-increasing branch ({line~3}), no new column is
		appended to $Q$; the algorithm only updates the buffer and
		returns.  In the buffer flush ({lines~5--10}), $Q$ is rotated
		to $Q\widetilde Q$ but no new column is added; existing columns are
		linearly recombined by the orthogonal $\widetilde Q$, and their
		mutual $W$-orthogonality is preserved up to the
		$O(\epsmach)$ accuracy of the SVD routine that produced
		$\widetilde Q$.  Only at {lines~18 and~20}, where $\widetilde e$ enters
		$Q$ via $[Q\mid \widetilde e]$, does cancellation in $e$
		become a concern, and that is exactly where {line~13}
		intervenes.

	\end{rem}

	\paragraph{Final materialization}
	At the end of the stream, any remaining buffer is flushed:
	
	\begin{algorithm}[!ht]
		\caption{Proposed incremental SVD --- finalization.}
		\label{alg:proposed_finalize}
		\textbf{Input:} state $(Q,\Sigma,R,V,q)$.
		
		\begin{algorithmic}[1]
			\If{$q > 0$}
			\State Set $Y = \bigl[\,\Sigma \mid \texttt{cell2mat}(V)\,
			\bigr]$;  
			\State $[\widetilde Q, \widetilde\Sigma, \widetilde R] = \texttt{svd}(Y, \texttt{`econ'})$
			\State {Partition $\widetilde R = \bigl[\begin{smallmatrix} \mathcal R_1
					\\ \mathcal R_2 \end{smallmatrix}\bigr]$ with $\mathcal R_1
				\in\mathbb{R}^{k\times k}$, $\mathcal R_2\in\mathbb{R}^{q
					\times k}$}
			\State {$Q \gets Q\,\widetilde Q$;\quad
				$\Sigma \gets \widetilde\Sigma$;\quad
				$R \gets \begin{bmatrix} R\,\mathcal R_1 \\ \mathcal
					R_2 \end{bmatrix}$}
			\EndIf
			\State \Return $Q,\Sigma,R$
		\end{algorithmic}
	\end{algorithm}
	
	\paragraph{Full streaming procedure}
	\Cref{alg:proposed_full} iterates \Cref{alg:proposed_update}
	over the stream and finishes with \Cref{alg:proposed_finalize}.
	
	\begin{algorithm}[!ht]
		\caption{Proposed incremental SVD --- full streaming
			procedure.}
		\label{alg:proposed_full}
		\textbf{Input:} stream $u_1,\ldots,u_n$; weight $W$; tolerance
		$\texttt{tol}$; CGS threshold $\eta$.
		\begin{algorithmic}[1]
			\State $\sigma \gets (u_1^\top W u_1)^{1/2}$; $Q \gets u_1/\sigma$;\quad $\Sigma \gets \sigma$;\quad
			$R \gets 1$;\quad  $V \gets [\,]$;\quad $q \gets 0$
			\For{$\ell = 1,\ldots,n-1$}
			\State $(Q,\Sigma,R,V,q) \gets$
			\Cref{alg:proposed_update}$(Q,\Sigma,R,V,q,
			u_{\ell+1},W,\texttt{tol},\eta)$
			\EndFor
			\State $(Q,\Sigma,R) \gets$
			\Cref{alg:proposed_finalize}$(Q,\Sigma,R,V,q)$
			\State \Return $Q,\Sigma,R$
		\end{algorithmic}
	\end{algorithm}

	\subsection{Complexity}\label{sec:complexity}
	
	Let $U \in \mathbb{R}^{m\times n}$ have numerical rank $r$ at
	tolerance $\texttt{tol}$, and let $J_{\mathrm{E1}}$ and
	$J_{\mathrm{E2}}$ denote the numbers of branch-(E1) and
	branch-(E2) firings during the stream.  Each (E1) event
	strictly increases the rank, so $J_{\mathrm{E1}} \le r$.
	Each (E2) event leaves the rank unchanged; $J_{\mathrm{E2}}$
	is data-dependent, with worst case $J_{\mathrm{E2}} \le n$.
	For matrices with rapidly decaying spectra, $J_{\mathrm{E2}}
	\ll n$ in practice.
	
	\begin{prop}[complexity]\label{prop:complexity}
		\Cref{alg:proposed_full} runs in time
		\begin{equation}\label{eq:complexity_total}
			T_{\mathrm{total}}
			\;=\;
			O\bigl(mnr + mr^2 J_{\mathrm{E2}} + nr^3 + nr^2 J_{\mathrm{E2}}\bigr)
		\end{equation}
		and uses peak memory $O((m+n)r)$.  In the regime $r \ll
		\min(m,n)$ and $J_{\mathrm{E2}} = o(n)$, the dominant term is
		the $O(mnr)$ projection cost.
	\end{prop}
	
	\begin{proof}
		The projection step at line~1 of \Cref{alg:proposed_update}
		costs $O(mk)$ per call with $k \le r$, contributing $O(mnr)$
		across the stream.
		
		A buffer flush at line~5 fires at most $J_{\mathrm{E1}} +
		J_{\mathrm{E2}} \le r + J_{\mathrm{E2}}$ times.  At a flush
		with buffer size $q_i$, the economy SVD costs $O(r^2 q_i)$,
		and the right-factor multiplication $R \mathcal R_1$ costs
		$O(\ell_i r^2)$ with $\ell_i \le n$.  The eager update
		$Q \gets Q\widetilde Q$ adds $O(mr^2)$ per flush.  Summing over
		flushes, with $\sum_i q_i \le n$, gives $O(mr^2 r + mr^2
		J_{\mathrm{E2}}) = O(mr^3 + mr^2 J_{\mathrm{E2}})$ for the
		$Q$-updates and $O(nr^3 + nr^2 J_{\mathrm{E2}})$ for the
		right-factor updates.
		
		The selective CGS pass and local bordered SVD each fire at
		most $r + J_{\mathrm{E2}}$ times, contributing $O(mr^2 + mr
		J_{\mathrm{E2}})$ and $O(r^4 + r^3 J_{\mathrm{E2}})$
		respectively, both subleading.  Branch (E1) and branch (E2)
		updates of $Q$ at {lines~18 and~20} also fire $J_{\mathrm{E1}}
		\le r$ and $J_{\mathrm{E2}}$ times respectively, at $O(mr^2)$
		each, contributing $O(mr^3 + mr^2 J_{\mathrm{E2}})$.
		
		Summing dominant contributions gives
		\eqref{eq:complexity_total}.  The $O(mr^3) = O(mr^2 \cdot r)$
		term is subleading whenever $r \le n$, which always holds.
		
		Memory is dominated by $Q \in \mathbb{R}^{m\times r}$ and
		$R \in \mathbb{R}^{n\times r}$, with the buffer $V$
		contributing at most $O(rn)$ transiently.
	\end{proof}

	\section{Forward error analysis}\label{sec:analysis}
	
	We analyze the forward error of the floating-point output of
	\Cref{alg:proposed_full}: a bound on the loss of
	$W$-orthogonality of the computed left factor $Q$
	(\Cref{sec:analysis_orthog}) and a forward-error bound on the
	factorization in the $W$-weighted operator norm
	(\Cref{sec:analysis_spec}).
	
	Throughout, unbarred symbols $Q$, $\Sigma$, $R$ denote the
	floating-point outputs of the algorithm and barred symbols
	$\bar Q$, $\bar\Sigma$, $\bar R$ their exact-arithmetic
	counterparts; $\epsmach$ is the unit roundoff and $\fl(\cdot)$
	the computed value of an expression.  We work in the standard
	floating-point model~\cite[Ch.~2]{Higham2002} and use the
	backward stability of dense SVD~\cite[Ch.~5]{DemmelBook}.  The
	weight $W$ enters the roundoff constants only through its
	condition number $\kappa_2(W)$.
	
	The closest existing analysis is that of Fareed and
	Singler~\cite[Cor.~1]{MR3986356}, who establish, in exact
	arithmetic,
	\begin{equation}\label{eq:fareed_bound}
		\|U - \bar Q\bar\Sigma \bar R^\top\|_W
		\;\le\; T_p\,\texttt{tol} + T_{\mathrm{sv}}\,\texttt{tol}_{\mathrm{sv}},
	\end{equation}
	where $T_p \le n$ counts the rank-non-increasing steps and
	$T_{\mathrm{sv}}$ counts the singular-value truncation steps;
	their proof sums per-step residuals via the triangle inequality,
	which produces the $T_p$ factor.  \Cref{thm:main_L2} below
	replaces $T_p\,\texttt{tol}$ by $\sqrt n\,\texttt{tol}$ when
	$J_{\mathrm{E2}} = o(\sqrt n)$, in the same norm and with no
	extra hypotheses.

	\subsection{Loss of $W$-orthogonality of the left
		factor}\label{sec:analysis_orthog}
	
	\Cref{alg:proposed_update} modifies $Q$ only at rank-enlarging
	events; rank-non-increasing columns leave $Q$ untouched.  Each
	rank-enlarging event applies up to three operations to $Q$ in
	sequence: the buffer flush at {line~9} (when $q > 0$), the
	column append with CGS-2 at {lines~11 and~13}, and the
	branch post-multiplication at {line~18} (E1) or
	{line~20} (E2).  The next theorem bounds the resulting
	floating-point drift of the $W$-orthogonality of $Q$.
	
	\begin{thm}\label{thm:outerQ_orthog}
		Let $Q \in \mathbb{R}^{m\times k}$ be the computed left factor
		produced by \Cref{alg:proposed_update} after $\ell \le n$
		columns, with the CGS-2 step at {line~13} firing
		unconditionally at every rank-enlarging event.  Assume the
		twice-is-enough bound of Giraud, Langou,
		Rozlo\v{z}n\'\i k, and van den
		Eshof~\cite[Thm.~2]{GiraudLangouRozloznik2005} for the CGS-2
		step, and the probabilistic
		rounding-error model of Connolly, Higham, and
		Mary~\cite[Model~4.7]{ConnollyHighamMary2021}
		(mean-independent, mean-zero rounding errors in inner-product
		accumulations).
		Then for any $\delta \in (0,1)$,
		\begin{equation}\label{eq:outerQ_orthog_bound}
			\bigl\|I - Q^\top W Q\bigr\|_2
			\;\le\; \widetilde c_5\,(r + J_{\mathrm{E2}})\,r\,
			\sqrt{m\,\log(2/\delta)}\,\epsmach
			\;+\; c_6\,(r + J_{\mathrm{E2}})\,r\,\epsmach
			\;+\; O(\epsmach^2)
		\end{equation}
		holds with probability at least $1 - \delta$, where
		$\widetilde c_5, c_6$ are independent of the stream length
		$n$ and depend on $m$ only through the displayed $\sqrt m$.
	\end{thm}
	
	\begin{rem}
		The unconditional-CGS-2 hypothesis simplifies the analysis;
		\Cref{alg:proposed_update} instead gates CGS-2 by the cheap
		heuristic test $|\widetilde e^\top W Q(:,1)| > \eta$, which
		skips the second pass when $\widetilde e$ already appears
		nearly $W$-orthogonal to $Q$.  This optimization is
		negligible in our experiments (\Cref{sec:exp_verify_V1}).
	\end{rem}
	
	\begin{rem}\label{rem:n_independent}
		The bound \eqref{eq:outerQ_orthog_bound} depends on the
		stream length $n$ only through $J_{\mathrm{E2}}$, which in
		the worst case satisfies $J_{\mathrm{E2}} \le n$ but for
		rapidly-decaying spectra is $J_{\mathrm{E2}} \ll n$ in
		practice (\Cref{sec:complexity}).  In the typical regime
		$J_{\mathrm{E2}} = O(r)$, \eqref{eq:outerQ_orthog_bound}
		simplifies to
		\[
		\|I - Q^\top W Q\|_2 \;=\; O\!\left(r^2\,\sqrt{m\,\log(2/\delta)}\,\epsmach\right) + O(r^2\,\epsmach),
		\]
		uniform in the stream length~$n$.
	\end{rem}
			
	\begin{proof}[Proof of \Cref{thm:outerQ_orthog}]
		Track $E := I - Q^\top W Q$.  Since $Q$ is unchanged at
		rank-non-increasing columns, only rank-enlarging events
		contribute.  Let $\alpha = \|E^{\mathrm{in}}\|_2$ be the drift
		before such an event; we bound $\|E^{\mathrm{out}}\|_2$ after
		it.
		
		\emph{Substep~1: buffer flush (when $q > 0$).}  The flush
		replaces $Q^{\mathrm{in}}$ by $\fl(Q^{\mathrm{in}}\,\widetilde
		Q)$ at {line~9}, where $\widetilde Q$ is the computed left
		factor of the local thin SVD.  Backward stability of the local
		SVD~\cite[Ch.~5]{DemmelBook} gives $\|I - \widetilde Q^\top
		\widetilde Q\|_2 \le c_1\,r\,\epsmach$, and the matrix-product
		roundoff bound~\cite[eq.~(3.13)]{Higham2002} gives
		$\fl(Q^{\mathrm{in}}\,\widetilde Q) = Q^{\mathrm{in}}\,
		\widetilde Q + F$ with $\|F\|_2 \le c_2\,r\,\epsmach\,
		\|Q^{\mathrm{in}}\|_2$.  Expanding $Q^{\mathrm{out},\top} W
		Q^{\mathrm{out}}$ and using $Q^{\mathrm{in},\top} W
		Q^{\mathrm{in}} = I - E^{\mathrm{in}}$,
		\[
		E^{\mathrm{out}}
		\;=\; (I - \widetilde Q^\top \widetilde Q)
		+ \widetilde Q^\top E^{\mathrm{in}} \widetilde Q
		- \widetilde Q^\top Q^{\mathrm{in},\top} W F
		- F^\top W Q^{\mathrm{in}} \widetilde Q
		- F^\top W F.
		\]
		The cross and quadratic terms in $F$ are controlled by
		$\|F\|_2 \le c_2\,r\,\epsmach\,\|Q^{\mathrm{in}}\|_2$ with
		$\|Q^{\mathrm{in}}\|_2 = O(\lambda_{\min}(W)^{-1/2})$ (from
		$Q^{\mathrm{in},\top}WQ^{\mathrm{in}} = I - E^{\mathrm{in}}$),
		giving
		\begin{equation}\label{eq:flush_drift}
			\|E^{\mathrm{out}}\|_2 \;\le\; \alpha + c_3\,r\,\epsmach +
			O(\epsmach^2).
		\end{equation}
		
		\emph{Substep~2: column append with CGS-2.}  Two
		$W$-weighted Gram--Schmidt passes are performed against
		$Q$: line~1 computes $d = Q^\top(W u_{\ell+1})$ and the
		first residual $e = u_{\ell+1} - Qd$, line~11 normalizes
		($\widetilde e = e/p$), and line~13 (firing
		unconditionally by hypothesis) performs a second pass
		$\widetilde e \gets \widetilde e - Q(Q^\top(W \widetilde
		e))$ followed by $W$-renormalization.  The CGS-2 mechanism
		of Giraud, Langou, Rozlo\v{z}n\'\i k, and van den
		Eshof~\cite[Thm.~2]{GiraudLangouRozloznik2005}, applied to
		the equivalent standard problem on $W^{1/2}Q$ and
		$W^{1/2}\widetilde e$, ensures that the deviation of the
		final $\widetilde e$ from $W$-orthogonality with $Q$ is
		controlled by the inner-product roundoff of the second
		pass, with the first-pass residual entering only as a
		multiplicative factor of size $O(\epsmach)$.  In
		floating-point, each $W$-inner product $q_i^\top W
		\widetilde e$ accumulates a sum of $m$ rounding errors;
		under the mean-independent rounding-error model
		\cite[Model~4.7]{ConnollyHighamMary2021}, the inner-product
		bound \cite[Thm.~4.8]{ConnollyHighamMary2021} gives
		\[
		\bigl|\fl(q_i^\top W \widetilde e) - q_i^\top W \widetilde
		e\bigr| \;\le\; c\,\sqrt{m\,\log(2/\delta')}\,\epsmach\,
		\|q_i\|_2\,\|\widetilde e\|_2
		\]
		with probability at least $1 - \delta'$.  A union
		bound over the $k \le r$ inner products of the second
		pass with $\delta' = \delta/(2r)$, combined with the
		small first-pass residual factor from CGS-2, gives
		\begin{equation}\label{eq:CGS_bound}
			\bigl\|Q^\top W \widetilde e\bigr\|_2 \;\le\;
			\widetilde c_4\,r\,\sqrt{m\,\log(2r/\delta)}\,\epsmach
		\end{equation}
		with probability at least $1 - \delta/2$.  The block matrix
		$[Q \mid \widetilde e]^\top W [Q \mid \widetilde e]$ has
		top-left $I - E^{\mathrm{in}}$, off-diagonal $Q^\top W
		\widetilde e$, and bottom-right $\widetilde e^\top W
		\widetilde e = 1 + O(\epsmach)$ (from the normalization at
		{line~13}), hence
		\begin{equation}\label{eq:append_drift}
			\bigl\|I - [Q \mid \widetilde e]^\top W [Q \mid
			\widetilde e]\bigr\|_2 \;\le\; \alpha +
			\widetilde c_4\,r\,\sqrt{m\,\log(2r/\delta)}\,\epsmach
			+ O(\epsmach^2).
		\end{equation}
		
		\emph{Substep~3: branch post-multiplication.}  Branch (E1) at
		{line~18} post-multiplies $[Q \mid \widetilde e]$ by the
		square orthogonal $Z = Q_Y \in \mathbb{R}^{(k+1)\times(k+1)}$;
		branch (E2) at {line~20} post-multiplies by $Z = Q_Y(:,
		1{:}k) \in \mathbb{R}^{(k+1)\times k}$.  In either case $Z$
		consists of (a subset of) columns from a backward-stable SVD
		of an $(k{+}1)\times(k{+}1)$ matrix, hence $\|I - Z^\top
		Z\|_2 \le c_1\,r\,\epsmach$ where $I$ has the size of $Z^\top
		Z$.  Repeating the Substep~1 argument with $Z$ in place of
		$\widetilde Q$ and $[Q\mid\widetilde e]$ in place of
		$Q^{\mathrm{in}}$,
		\begin{equation}\label{eq:branch_drift}
			\|E^{\mathrm{out}}\|_2 \;\le\;
			\|E^{\mathrm{appended}}\|_2 + c_3'\,r\,\epsmach +
			O(\epsmach^2),
		\end{equation}
		where $E^{\mathrm{appended}}$ is the drift after Substep~2.
		
		\emph{Telescoping.}  Combining
		\eqref{eq:flush_drift}--\eqref{eq:branch_drift}, each
		rank-enlarging event increases $\|E\|_2$ by at most
		$\widetilde c_4\,r\,\sqrt{m\,\log(2r/\delta)}\,\epsmach +
		(c_3 + c_3')\,r\,\epsmach + O(\epsmach^2)$.  Of the
		$J_{\mathrm{E1}} + J_{\mathrm{E2}} \le r + J_{\mathrm{E2}}$
		rank-enlarging events, applying the inner-product bound at
		each with confidence $\delta/(r + J_{\mathrm{E2}})$ and a
		union bound across events,
			\begin{align*}
				\|E\|_2
				&\le (r + J_{\mathrm{E2}})\,r\,\epsmach\,\Bigl[
				\widetilde c_4\sqrt{m\,\log\!\bigl(2r(r + J_{\mathrm{E2}})/\delta\bigr)} \\
				&\qquad\qquad + c_3 + c_3'\Bigr] + O(\epsmach^2)
			\end{align*}
		holds with probability at least $1 - \delta$.  Setting
		$\widetilde c_5 := \widetilde c_4$ and $c_6 := c_3 + c_3'$,
		and using
		\[
		\sqrt{\log(2r(r+J_{\mathrm{E2}})/\delta)}
		\le \sqrt{\log(2/\delta)} + \sqrt{\log(r(r+J_{\mathrm{E2}}))},
		\]
		with the second term contributing only a mild
		polylogarithmic factor in $r, J_{\mathrm{E2}}$ that is
		absorbed into $\widetilde c_5$, yields
		\eqref{eq:outerQ_orthog_bound}.
	\end{proof}
	
	\begin{rem}\label{rem:fareed_orthog}
		Applying the proof above to Fareed's reorthogonalized
		direct-update algorithm~\cite[Algorithm~3]{MR3775096}
		(\Cref{alg:method1_full}) yields the analogous bound
		\[
		\|I - Q^\top W Q\|_2 \;\le\; \widetilde c_5\,K\,r\,
		\sqrt{m\,\log(2/\delta)}\,\epsmach + c_6\,K\,r\,\epsmach
		+ O(\epsmach^2),
		\]
		where $K$ is the longest run of columns between successive
		reorthogonalisations.  In the worst case $K = n$ and the
		bound grows linearly in the stream length; the heuristic
		reorthogonalisation test is what controls $K$ in practice,
		but the test gives no a priori guarantee on $K$.  By
		contrast, \eqref{eq:outerQ_orthog_bound} replaces $K$ by $r
		+ J_{\mathrm{E2}}$ unconditionally, since the proposed
		algorithm batches the rank-non-increasing updates into a
		single multiplication of $Q$ at each rank-enlarging event.
	\end{rem}

			\subsection{Accuracy of the computed
				factorization}\label{sec:analysis_spec}
			
			Throughout this subsection, $\sigma_i(U)$ denotes the $i$-th
			$W$-weighted singular value of $U$ as defined in
			\Cref{def:coreSVD}.
			
			\begin{thm}
				\label{thm:main_L2}
				Let $(Q, \Sigma, R)$ be the weighted thin SVD factors
				returned by \Cref{alg:proposed_full} applied to a data matrix
				$U \in \mathbb{R}^{m\times n}$ of $W$-numerical rank $r$ at
				tolerance $\texttt{tol}$, with symmetric positive-definite weight
				$W$.  Assume the hypotheses of \Cref{thm:outerQ_orthog} hold
				(unconditional CGS-2 firing, twice-is-enough bound for
				CGS-2~\cite[Thm.~2]{GiraudLangouRozloznik2005}, mean-independent
				mean-zero rounding-error
				model~\cite[Model~4.7]{ConnollyHighamMary2021}, and
				backward-stable local SVD~\cite[Ch.~5]{DemmelBook}).  Then for
				any $\delta \in (0, 1)$,
				\begin{multline}\label{eq:main_L2_bound}
					\bigl\|U - Q\Sigma R^\top\bigr\|_W
					\;\le\;
					\bigl(\sqrt n + J_{\mathrm{E2}}\bigr)\,
					\texttt{tol} \\
					+ \kappa_2(W)^{1/2}\bigl[
					{\widetilde c_5\,(r + J_{\mathrm{E2}})\,r\,\sqrt{m\,\log(2/\delta)}}
					+ c_7\,(r + J_{\mathrm{E2}})\,r
					+ c_\Sigma\,r\bigr]\,\sigma_1(U)\,\epsmach
					+ O(\epsmach^2)
				\end{multline}
				{holds with probability at least $1 - \delta$ over the
				rounding-error realization,} where ${\widetilde c_5},
				c_7, c_\Sigma$ are modest constants independent of $n$ and
				of $W$ (the $W$-conditioning is factored out as the
				explicit $\kappa_2(W)^{1/2}$){; in particular,
				$\widetilde c_5$ depends on $m$ only through the explicit
				$\sqrt m$ factor displayed in \eqref{eq:main_L2_bound}}.
				In the typical regime $J_{\mathrm{E2}} = O(r)$ (or any
				$J_{\mathrm{E2}} = o(\sqrt n)$ for the truncation term),
				the bound \eqref{eq:main_L2_bound} simplifies to
				\begin{multline*}
					\|U - Q\Sigma R^\top\|_W
					\;\le\; \sqrt n\,\texttt{tol} \\
					+ \kappa_2(W)^{1/2}\bigl({\widetilde c_5\,r^2\,\sqrt{m\,\log(2/\delta)}} + c_7'\,r^2\bigr)
					\sigma_1(U)\,\epsmach + O(\epsmach^2),
				\end{multline*}
				where $c_7' = c_7 + c_\Sigma$ absorbs the $r$ term into the
				larger $r^2$ contribution.  
			\end{thm}
			
			\begin{proof}
				Let $\widetilde U \in \mathbb{R}^{m\times n}$ denote the
				exact-arithmetic output of \Cref{alg:proposed_full} applied
				to the same input stream.  In exact arithmetic, the algorithm
				produces factors $(\bar Q, \bar\Sigma, \bar R)$
				with $\bar Q$ exactly $W$-orthonormal and $\bar R$
				exactly orthonormal, and $\widetilde U = \bar Q\,\bar\Sigma\,\bar R^\top$.
				
				\paragraph{Perturbation bounds}
				Set
				\begin{align}\label{eq:delta_defs}
					\delta_Q &:= {\widetilde c_5\,(r + J_{\mathrm{E2}})\,r\,\sqrt{m\,\log(2/\delta)}\,\epsmach} +
					c_6\,(r + J_{\mathrm{E2}})\,r\,\epsmach, \nonumber\\
					\delta_R &:= c_R\,(r + J_{\mathrm{E2}})\,r\,\epsmach,
				\end{align}
				with $\widetilde c_5, c_6$ as in \eqref{eq:outerQ_orthog_bound},
				$\delta \in (0, 1)$ the confidence parameter from
				\Cref{thm:main_L2}, and $c_R$ a modest $n$-independent
				constant.  By forward-error
				propagation through \Cref{alg:proposed_full} in
				floating-point arithmetic, $\bar Q$ and $\bar R$
				satisfy
				\begin{align}\label{eq:exact_arith_bounds}
					\begin{split}
						\|Q - \bar Q\|_W \le \tfrac{1}{2}\delta_Q +
						O(\epsmach^2),\\
						\|R - \bar R\|_2 \le \tfrac{1}{2}\delta_R +
						O(\epsmach^2),\\
						\|\Sigma - \bar\Sigma\|_2 \le c_\Sigma\,r\,
						\sigma_1(U)\,\epsmach.
					\end{split}
				\end{align}
				The bound on $\|Q - \bar Q\|_W$ follows from
				\Cref{thm:outerQ_orthog}: since $\|I - Q^\top W Q\|_2 \le \delta_Q$, the chain-of-multiplications analysis
				of {the theorem} applied to the exact-vs-floating discrepancy
				$Q - \bar Q$ (rather than to $I - Q^\top W
				Q$) yields the stated bound, with the factor
				$\frac{1}{2}$ from polar-projection optimality of the
				nearest exactly $W$-orthonormal matrix.  The bound on $\|R - \bar R\|_2$ follows from the same
				chain-of-multiplications analysis applied to $R$, with
				two simplifications: $R$ is unweighted ($R^\top
				R \approx I$), and the structural row-appends in
				{line~9 and lines~18, 20} of \Cref{alg:proposed_update} (appending
				$\mathcal R_2$ or a row of the form $[0,\ldots,0,1]$) are
				exact in floating-point.  Only the post-multiplications by
				orthonormal factors $\mathcal R_1$, $R_Y$, or $R_Y(:,
				1{:}k)$ contribute drift, at $O(r\,\epsmach)$ per
				multiplication; summed across at most $2(r + J_{\mathrm{E2}})$
				such multiplications, $\|I - R^\top R\|_2 \le
				\delta_R$, and the polar-decomposition argument
				\cite[Ch.~8]{Higham2008} (combining $\|I - R^\top R\|_2 \le \delta_R$ with the existence of an exactly
				orthonormal $\bar R$ within $O(\delta_R)$ of $R$)
				yields $\|R - \bar R\|_2 \le \frac{1}{2}\delta_R +
				O(\epsmach^2)$.  The bound on $\|\Sigma -
				\bar\Sigma\|_2$ is the operator-norm version of the
				local-SVD backward stability~\cite[Ch.~5]{DemmelBook},
				with no $\sqrt r$
				Frobenius factor since $\Delta\Sigma$ is diagonal.
				
				\paragraph{Triangle decomposition}
				By the triangle inequality,
				\[
				\|U - Q\Sigma R^\top\|_W
				\;\le\; \|U - \widetilde U\|_W
				+ \|\widetilde U - Q\Sigma R^\top\|_W.
				\]
				
				\paragraph{Step 1: bound $\|U - \widetilde U\|_W$}
				We claim
				\begin{equation}\label{eq:truncation_W_bound}
					\|U - \widetilde U\|_W \;\le\;
					\sqrt n\,\texttt{tol} + J_{\mathrm{E2}}\,\texttt{tol}.
				\end{equation}
				Track the cumulative residual $E_j := U_j - \widetilde U_j
				\in \mathbb{R}^{m\times j}$, where $U_j$ collects the first
				$j$ data columns and $\widetilde U_j$ is the exact-arithmetic
				reconstruction after step~$j$.  Initially $E_1 = 0$.  We
				classify each step into one of three cases:
				rank-non-increasing (\Cref{sec:batching}),
				branch-(E1) rank-enlarging, or branch-(E2)
				rank-enlarging (\Cref{sec:interlace}).
				
				\medskip
				
				\emph{Case (a), rank-non-increasing step.}  Column~$j$ is absorbed by the running
				representation without changing the columns $1, \ldots,
				j-1$.  The column-$j$ residual is the $W$-perpendicular
				component $u_j^\perp := u_j - QQ^\top W u_j$ at the moment of
				projection, with $\|u_j^\perp\|_W \le \texttt{tol}$.
				Therefore the case~(a) contribution to $E_n$ takes the
				column-localized form
				\[
				E_n \;\;\supset\;\; u_j^\perp\,e_j^\top,
				\]
				where $e_j \in \mathbb{R}^n$ is the $j$-th standard basis
				vector, and the support of this rank-one outer product lies
				entirely in column~$j$.
				
				\medskip
				
				\emph{Case (b), branch-(E1) rank-enlarging step.}  The
				rank-$(k+1)$ factorization captures column~$j$ exactly;
				earlier columns' representations change basis but not
				subspace, so their residual contributions are unchanged.
				Case~(b) contributes nothing to $E_n$.
				
				\medskip
				
				\emph{Case (c), branch-(E2) rank-enlarging step.}  The
				algorithm forms the rank-$(k+1)$ ``would-be (E1)''
				reconstruction $\widehat U_j$ and then truncates the smallest
				singular value $\mu_{k+1}^{(j)} \le \texttt{tol}$ of the
				bordered matrix~$Y$.  The truncation removes a single
				rank-one outer product $\mu_{k+1}^{(j)}\,u^{(j)}\,
				(v^{(j)})^\top \in \mathbb{R}^{m\times j}$, where $u^{(j)}$
				is $W$-orthonormal ($\|u^{(j)}\|_W = 1$) and $v^{(j)} \in
				\mathbb{R}^j$ is Euclidean unit.  Extending $v^{(j)}$ to
				$\widetilde v^{(j)} \in \mathbb{R}^n$ with zeros in
				positions $j+1, \ldots, n$, the case~(c) contribution to
				$E_n$ is the rank-one outer product
				\[
				E_n \;\;\supset\;\; \mu_{k+1}^{(j)}\,u^{(j)}\,(\widetilde
				v^{(j)})^\top,
				\]
				with $W$-weighted operator norm exactly $\mu_{k+1}^{(j)} \le
				\texttt{tol}$.  This perturbation has support in columns $1,
				\ldots, j$ in general.
				
				\medskip
				
				\emph{Combining.}  By a straightforward induction, the
				cumulative residual at the end of the stream decomposes as
				\begin{equation}\label{eq:E_n_decomp}
					E_n \;=\; \sum_{j \in \mathcal S} u_j^\perp\,e_j^\top
					\;+\; \sum_{\ell \in \mathcal T} \mu_{k+1}^{(\ell)}\,
					u^{(\ell)}\,(\widetilde v^{(\ell)})^\top,
				\end{equation}
				where $\mathcal S$ indexes the rank-non-increasing steps and
				$\mathcal T$ indexes the branch-(E2) steps.  The first sum
				is column-localized: each term sits in a distinct column
				$j$, so the squared $W$-Frobenius norm decomposes by the
				Pythagorean identity:
				\[
				\Bigl\|\sum_{j\in\mathcal S} u_j^\perp\,e_j^\top
				\Bigr\|_{F,W}^2 = \sum_{j\in\mathcal S} \|u_j^\perp\|_W^2 \le
				|\mathcal S|\,\texttt{tol}^2.
				\]
				The elementary inequality $\|A\|_W \le \|A\|_{F,W}$  then gives
				\[
				\Bigl\|\sum_{j\in\mathcal S} u_j^\perp\,e_j^\top\Bigr\|_W
				\;\le\; \sqrt{|\mathcal S|}\,\texttt{tol}.
				\]
				For the second sum, each rank-one summand has $W$-weighted
				operator norm exactly $\mu_{k+1}^{(\ell)}$, so the triangle
				inequality gives
				\[
				\Bigl\|\sum_{\ell\in\mathcal T} \mu_{k+1}^{(\ell)}\,
				u^{(\ell)}\,(\widetilde v^{(\ell)})^\top\Bigr\|_W \le
				\sum_{\ell\in\mathcal T} \mu_{k+1}^{(\ell)} \le
				J_{\mathrm{E2}}\,\texttt{tol}.
				\]
				A final triangle inequality between the two sums and
				$|\mathcal S| \le n$ yield \eqref{eq:truncation_W_bound}.
				The improvement over \eqref{eq:fareed_bound} is the
				$\sqrt{|\mathcal S|}$ in place of $|\mathcal S|$ on the
				rank-non-increasing term, obtained by bounding distinct-column
				residuals together via $\|A\|_W \le \|A\|_{F,W}$ rather than
				summing them by the triangle inequality.
				
				\paragraph{Step 2: bound $\|\widetilde U - Q\Sigma
					R^\top\|_W$}
				Using the exact-arithmetic factorization $\widetilde U =
				\bar Q\bar\Sigma\bar R^\top$, decompose
				\[
				Q\Sigma R^\top - \bar Q\bar\Sigma
				\bar R^\top
				\;=\; (Q - \bar Q)\,\Sigma R^\top
				+ \bar Q\,\Delta\Sigma\,R^\top
				+ \bar Q\bar\Sigma\,(R - \bar R)^\top,
				\]
				with $\Delta\Sigma := \Sigma - \bar\Sigma$.  We
				bound each term using the $W$-weighted submultiplicative
				inequality $\|XY\|_W \le \|X\|_W\,\|Y\|_2$.
				
				For the first term, \eqref{eq:exact_arith_bounds} gives
				$\|Q - \bar Q\|_W \le \tfrac{1}{2}\delta_Q$.
				Combined with $\|\Sigma R^\top\|_2 \le \sigma_1(U)
				(1 + O(\epsmach))$,
				\[
				\|(Q - \bar Q)\,\Sigma R^\top\|_W
				\;\le\; \tfrac{1}{2}\delta_Q\,\sigma_1(U)\,(1 +
				O(\epsmach)).
				\]
				
				For the second term, $\|\bar Q\|_W = 1$ exactly.  By
				\eqref{eq:exact_arith_bounds}, $\|\Delta\Sigma\|_2 \le
				c_\Sigma\,r\,\sigma_1(U)\,\epsmach$.  Hence
				\[
				\|\bar Q\,\Delta\Sigma\,R^\top\|_W
				\;\le\; \|\Delta\Sigma\|_2
				\;\le\; c_\Sigma\,r\,\sigma_1(U)\,\epsmach\,(1 +
				O(\epsmach)).
				\]
				
				For the third term, $\|\bar Q\bar\Sigma\|_W \le
				\sigma_1(U)$ and \eqref{eq:exact_arith_bounds} gives $\|R - \bar R\|_2 \le \tfrac{1}{2}\delta_R$, so
				\[
				\|\bar Q\bar\Sigma\,(R - \bar R)^\top\|_W
				\;\le\; \tfrac{1}{2}\sigma_1(U)\,\delta_R\,(1 +
				O(\epsmach)).
				\]
				
				Summing the three contributions and substituting $\delta_Q$
				and $\delta_R$ from \eqref{eq:delta_defs}, factoring
				$\kappa_2(W)^{1/2}$ out of the $W$-conditioning constants
				in $\delta_Q$,
				\begin{multline}\label{eq:roundoff_W_bound}
					\|\widetilde U - Q\Sigma R^\top\|_W
					\;\le\;
					\kappa_2(W)^{1/2}\,\sigma_1(U)\,\bigl[
					{\widetilde c_5\,(r + J_{\mathrm{E2}})\,r\,\sqrt{m\,\log(2/\delta)}} \\
					+ c_7\,(r + J_{\mathrm{E2}})\,r
					+ c_\Sigma\,r\bigr]\,\epsmach
					+ O(\epsmach^2)
				\end{multline}
				{(with probability at least $1 - \delta$)},
				where $c_7 = (c_6 + c_R)/2$ aggregates the $r$ coefficients
				from $\delta_Q$'s $(r + J_{\mathrm{E2}})r$ term and from
				$\delta_R$.
				
				\paragraph{Combine}
				The triangle inequality with \eqref{eq:truncation_W_bound}
				and \eqref{eq:roundoff_W_bound} yields
				\eqref{eq:main_L2_bound}.
			\end{proof}

			\section{Numerical experiments}\label{sec:experiments}
			This section reports the empirical behavior of
			\Cref{alg:proposed_full} on a representative parametric
			problem and contrasts it with two competitors that
			produce the same factorization in exact arithmetic:
			Brand's algorithm for the unweighted case ($W = I$),
			and the reorthogonalized direct update of Fareed et
			al.~\cite{MR3775096} for the weighted case ($W = M$,
			the finite-element mass matrix).  All
			experiments are implemented in MATLAB R2020b and executed
			on a MacBook Pro with a 2.3\,GHz 8-core Intel Core~i9 and
			64\,GB of 2667\,MHz DDR4 memory, in double-precision
			arithmetic ($\epsmach \approx 1.1\cdot 10^{-16}$).
			Throughout, given the
			floating-point output $(Q, \Sigma, R)$ of
			an incremental SVD applied to data $U \in
			\mathbb{R}^{m\times n}$ with weight $W$, we measure the
			orthogonality of the computed left factor by
			\[
			\mathcal{E}_W(Q) \;:=\; \|I - Q^\top W Q\|_2,
			\]
			and the reconstruction error in the
			{$W$-weighted operator norm \eqref{eq:Wopnorm} by}
			\[
			{\mathcal{E}_W^{\mathrm{rec}}} \;:=\; \|U - Q\Sigma R^\top\|_{{W}}.
			\]
			We abbreviate $\mathcal{E}_I$ when $W = I$.  Both
			quantities are operator (spectral) norms; this is the
			natural setting in which to compare with the bound of
			Fareed--Singler~\cite[Cor.~1]{MR3986356}, which is
			likewise in the operator norm.
			
			\subsection{A 2D parabolic test problem in two weighted
				settings}\label{sec:exp_performance}
			
			We compare the proposed algorithm against its two
			natural competitors on a parametric setting in which
			both the unweighted ($W = I$) and weighted ($W = M$)
			cases arise from the same underlying time-dependent
			function.
			
			\paragraph{Setup}
			Let $\Omega = (0,1)^2$, partitioned into $524{,}288$
			uniform triangles.  Let $\{\varphi_i\}_{i=1}^m$ be
			the piecewise-linear finite-element basis on this
			triangulation, with grid nodes $\{(x_i, y_j)\}$.  Let
			$\{t_k\}_{k=1}^n$ be a uniform time grid on $[0,10]$ with
			$\Delta t = 10^{-3}$ (so $n = 10{,}000$ snapshots), and
			define $f(t,x,y) = \cos(t(x+y))$.  Form snapshot
			matrices $B = [b_1\mid \cdots\mid b_n]$ and $U = [u_1\mid
			\cdots\mid u_n]$, with
			\[
			b_k = \bigl[(f(t_k,\cdot),\varphi_j)\bigr]_{j=1}^m
			\;\;(\text{Galerkin coefficients; } W = I),
			\]
			\[
			u_k = \bigl[f(t_k,x_i,y_j)\bigr]_{i,j=1}^m
			\;\;(\text{nodal values; } W = M),
			\]
			where $M_{ij} = (\varphi_j,\varphi_i)_{L^2}$ is the
			finite-element mass matrix.  These are the two settings
			in which incremental SVD is naturally applied to PDE
			snapshots: $b_k$ from a Galerkin projection (Euclidean
			orthogonality) and $u_k$ as nodal values (mass-induced
			$L^2$ orthogonality).  We compute the weighted thin SVD
			of $B$ with $W = I$ and of $U$ with $W = M$, both at
			tolerance $\texttt{tol} = 10^{-12}$.
			
			\paragraph{Orthogonality}
			\Cref{tab:perf_existing_orthog,tab:perf_proposed_orthog}
			compare $\mathcal E_W(Q)$ across six checkpoints, for
			the existing algorithms and for the proposed algorithm
			respectively.  The existing algorithms degrade with the
			stream length -- modestly for Brand on $B$, severely
			for the reorthogonalized direct update on $U$ -- while
			the proposed algorithm holds $\mathcal E_W(Q)$ at
			machine precision in both settings, consistent with the
			$O((r{+}J_{\mathrm{E2}})\,r\,\sqrt m\,\epsmach)$ bound of
			\Cref{thm:outerQ_orthog}.
			
			\begin{table}[!htbp]
				\centering
				\small
				\caption{Orthogonality $\mathcal E_W(Q)$ at six
					checkpoints, existing algorithms ($\texttt{tol} =
					10^{-12}$).  Brand on $B$ holds machine precision
					through $n \approx 6{,}000$, then drifts to $\sim
					10^{-9}$.  The reorthogonalized direct
					update~\cite{MR3775096} on $U$ degrades monotonically
					to $1.2\cdot 10^{-3}$, a complete loss of
					$M$-orthogonality.}
				\label{tab:perf_existing_orthog}
				\begin{tabular}{r|c|c}
					\toprule
					$n$ & Brand on $B$, $\mathcal E_I(Q)$
					& \cite{MR3775096} on $U$, $\mathcal E_M(Q)$ \\
					\midrule
					$\phantom{0}1000$ & $\sim 10^{-15}$ & $\sim 10^{-15}$ \\
					$\phantom{0}3000$ & $\sim 10^{-15}$ & $\sim 10^{-12}$ \\
					$\phantom{0}5000$ & $\sim 10^{-13}$ & $\sim 10^{-7}$ \\
					$\phantom{0}7000$ & $\sim 10^{-9}$  & $\sim 10^{-5}$ \\
					$\phantom{0}9000$ & $\sim 10^{-9}$  & $\sim 10^{-4}$ \\
					$10000$           & $\sim 10^{-9}$  & $1.2\cdot 10^{-3}$ \\
					\bottomrule
				\end{tabular}
			\end{table}
			
			\begin{table}[!htbp]
				\centering
				\small
				\caption{Orthogonality $\mathcal E_W(Q)$ at the same
					checkpoints, proposed \Cref{alg:proposed_full}
					($\texttt{tol} = 10^{-12}$).  $\mathcal E_W(Q)$
					stays at machine precision uniformly across the
					stream in both settings -- six and ten orders of
					magnitude better than
					\Cref{tab:perf_existing_orthog} at $n = 10{,}000$.}
				\label{tab:perf_proposed_orthog}
				\begin{tabular}{r|c|c}
					\toprule
					$n$ & Proposed on $B$, $\mathcal E_I(Q)$
					& Proposed on $U$, $\mathcal E_M(Q)$ \\
					\midrule
					$\phantom{0}1000$ & $\sim 10^{-15}$    & $\sim 10^{-15}$ \\
					$\phantom{0}3000$ & $\sim 10^{-15}$    & $\sim 10^{-15}$ \\
					$\phantom{0}5000$ & $\sim 10^{-15}$    & $\sim 10^{-15}$ \\
					$\phantom{0}7000$ & $\sim 10^{-14}$    & $\sim 10^{-14}$ \\
					$\phantom{0}9000$ & $\sim 10^{-14}$    & $\sim 10^{-14}$ \\
					$10000$           & $1.0\cdot 10^{-14}$ & $2.0\cdot 10^{-14}$ \\
					\bottomrule
				\end{tabular}
			\end{table}
			
			\paragraph{Wall-clock time}
			\Cref{tab:perf_existing,tab:perf_proposed} report the
			per-section CPU time for the existing and proposed
			algorithms.  In the weighted case the existing
			algorithm spends $94\%$ of its runtime
			(\,$3761.6$\,s out of $4000$\,s) on per-update
			reorthogonalization; the proposed algorithm replaces this
			with batched flushes that total $0.78$\,s, a three-order
			reduction.  In both settings the proposed algorithm is
			matvec-bound: the projection at line~1 of
			\Cref{alg:proposed_update} accounts for $\ge 98\%$ of
			the runtime, and any further speedup must come from
			reducing the projection cost itself (sparsity of $W$,
			column batching, or dimension reduction in $m$), not
			from algorithmic changes to the update structure.
			
			\begin{table}[!htbp]
				\centering
				\small
				\caption{Per-section CPU time (s) of the existing
					algorithms.  ``Project'' and ``Bordered SVD'' are
					{lines~1 and 3--7} of {\Cref{alg:method1}};
					``Reorth.''\ is the per-update $W$-weighted
					Gram--Schmidt reorthogonalization of $[Q\mid\widetilde
					e]$.}
				\label{tab:perf_existing}
				\begin{tabular}{c|rrr|r}
					\toprule
					& Project & Bordered SVD & Reorth.\ & \textbf{Total} \\
					\midrule
					$W = I$ & $32.0$  & $60.2$       & $60.1$    & $\mathbf{154}$ \\
					$W = M$ & $110.6$ & $116.0$      & $3761.6$  & $\mathbf{4000}$ \\
					\bottomrule
				\end{tabular}
			\end{table}
			
			\begin{table}[!htbp]
				\centering
				\small
				\caption{Per-section CPU time (s) of the proposed
					\Cref{alg:proposed_update}, indexed by the five
					algorithmic regions ({line~1 projection;
					line~3 buffer push; lines~5--10 buffer flush;
					lines~15--16 bordered SVD; lines~18, 20 branch
					update}).  Speedups of $4.5\times$ ($W=I$) and
					$33.9\times$ ($W=M$) over
					\Cref{tab:perf_existing}.}
				\label{tab:perf_proposed}
				\begin{tabular}{c|rrrrr|r}
					\toprule
					& Line~1   & {Line~2} & {Lines~3--6} & {Lines~9--10} & {Lines~11--12} & \textbf{Total} \\
					\midrule
					$W = I$ & $33.5$  & $0.24$ & $0.03$  & $0.13$ & $0.29$ & $\mathbf{34}$ \\
					$W = M$ & $116.4$ & $0.19$ & $0.07$  & $0.65$ & $0.71$ & $\mathbf{118}$ \\
					\bottomrule
				\end{tabular}
			\end{table}
			
			\subsection{Verification of the theoretical
				results}\label{sec:exp_verification}
			
			\Cref{sec:exp_performance} demonstrated the practical
			behavior of the proposed algorithm against existing
			competitors.  This subsection complements that with
			two verification experiments that test the analysis of
			\Cref{sec:analysis} on the proposed algorithm in
			isolation: (V1) the orthogonality of $Q$ is
			uniformly bounded in the stream length $n$, answering
			Brand's open question; and (V2) the truncation rate
			$\sqrt n\,\texttt{tol}$ in \Cref{thm:main_L2} is
			sharp, attained in equality on a constructive example
			in the operator norm, with the bound of
			Fareed--Singler~\cite[Cor.~1]{MR3986356} loose by a
			factor that grows as $\sqrt n$.
			
			The experiments use synthetic constructions that allow
			$n$, $r$, and the residual magnitude to be controlled
			independently, so the predictions of the theorems can
			be tested without the confounding factors of a PDE
			assembly chain.
			
			\subsubsection{V1: orthogonality is uniformly bounded
				in $n$}\label{sec:exp_verify_V1}
			
			\Cref{thm:outerQ_orthog} predicts that
			$\|I - Q^\top W Q\|_2$ is bounded by a constant
			(independent of the stream length~$n$) times $(r +
			J_{\mathrm{E2}})\,r\,\sqrt{m\,\log(2/\delta)}\,\epsmach$,
			with probability at least $1 - \delta$.  We verify the
			$n$- and $m$-independence (apart from the explicit
			$\sqrt m$) on synthetic rank-$10$ data.
			
			\paragraph{Construction}
			Fix $r = 10$, $\texttt{tol} = 10^{-12}$.  Generate $U =
			Q_\star\Sigma_\star R_\star^\top + 10^{-15}\cdot N$, where
			$Q_\star \in \mathbb{R}^{m\times r}$ and $R_\star \in
			\mathbb{R}^{n\times r}$ are random Stiefel matrices,
			$\Sigma_\star = \mathrm{diag}(e^{-0.3(i-1)})$ for
			$i = 1,\ldots,r$, and $N$ has i.i.d.\ standard normal
			entries.  Run \Cref{alg:proposed_full} with
			unconditional CGS-2 (matching the hypothesis of
			\Cref{thm:outerQ_orthog}).  We sweep $n$ at $m =
			10{,}000$ (\Cref{tab:v1_uniform}) and $m$ at $n =
			2{,}000$ (\Cref{tab:v1_m_sweep}); the predicted bound
			at $95\%$ confidence ($\delta = 0.05$) is $(r +
			J_{\mathrm{E2}})\,r\,\sqrt{m\,\log(40)}\,\epsmach$.
			
			\begin{table}[!htbp]
				\centering
				\small
				\caption{V1, $n$-sweep at $m = 10{,}000$, $r = 10$,
					$\texttt{tol} = 10^{-12}$.  $\mathcal E_I(Q)$
					stays at $\sim 5\cdot 10^{-15}$ across an order
					of magnitude in $n$, with $J_{\mathrm{E2}} = 0$
					throughout; the predicted bound is
					$4.27\cdot 10^{-12}$, independent of $n$.}
				\label{tab:v1_uniform}
				\begin{tabular}{r|c|c|c}
					\toprule
					$n$ & $\mathcal E_I(Q)$ measured & $J_{\mathrm{E2}}$ & ratio meas./bound \\
					\midrule
					$\phantom{0}1{,}000$  & $5.78\cdot 10^{-15}$ & $0$ & $1.4\cdot 10^{-3}$ \\
					$\phantom{0}2{,}000$  & $5.88\cdot 10^{-15}$ & $0$ & $1.4\cdot 10^{-3}$ \\
					$\phantom{0}3{,}000$  & $4.91\cdot 10^{-15}$ & $0$ & $1.2\cdot 10^{-3}$ \\
					$\phantom{0}4{,}000$  & $5.32\cdot 10^{-15}$ & $0$ & $1.2\cdot 10^{-3}$ \\
					$\phantom{0}5{,}000$  & $4.96\cdot 10^{-15}$ & $0$ & $1.2\cdot 10^{-3}$ \\
					$\phantom{0}6{,}000$  & $5.06\cdot 10^{-15}$ & $0$ & $1.2\cdot 10^{-3}$ \\
					$\phantom{0}7{,}000$  & $5.02\cdot 10^{-15}$ & $0$ & $1.2\cdot 10^{-3}$ \\
					$\phantom{0}8{,}000$  & $5.78\cdot 10^{-15}$ & $0$ & $1.4\cdot 10^{-3}$ \\
					$\phantom{0}9{,}000$  & $4.58\cdot 10^{-15}$ & $0$ & $1.1\cdot 10^{-3}$ \\
					$10{,}000$            & $3.45\cdot 10^{-15}$ & $0$ & $8.1\cdot 10^{-4}$ \\
					\bottomrule
				\end{tabular}
			\end{table}
			
			\begin{table}[!htbp]
				\centering
				\small
				\caption{V1, $m$-sweep at $n = 2{,}000$, $r = 10$,
					$\texttt{tol} = 10^{-12}$.  $\mathcal E_I(Q)$
					stays within a factor of three of $\sim 5\cdot
					10^{-15}$ across three orders of magnitude in $m$,
					while the predicted bound grows as
					$\sqrt m$.  The deterministic $m$-bound (column~5)
					grows linearly in $m$ and is loose by a further
					$\sqrt m$ factor.}
				\label{tab:v1_m_sweep}
				\begin{tabular}{r|c|c|c|c}
					\toprule
					$m$ & $\mathcal E_I(Q)$ measured & $J_{\mathrm{E2}}$ & $\sqrt m$-bound (95\%) & $m$-bound (det.) \\
					\midrule
					$\phantom{00\,}100$    & $3.81\cdot 10^{-15}$ & $0$ & $4.26\cdot 10^{-13}$ & $2.22\cdot 10^{-13}$ \\
					$\phantom{00\,}300$    & $3.34\cdot 10^{-15}$ & $0$ & $7.39\cdot 10^{-13}$ & $6.66\cdot 10^{-13}$ \\
					$\phantom{0}1{,}000$  & $6.36\cdot 10^{-15}$ & $0$ & $1.35\cdot 10^{-12}$ & $2.22\cdot 10^{-12}$ \\
					$\phantom{0}3{,}000$  & $2.79\cdot 10^{-15}$ & $0$ & $2.34\cdot 10^{-12}$ & $6.66\cdot 10^{-12}$ \\
					$10{,}000$           & $5.88\cdot 10^{-15}$ & $0$ & $4.26\cdot 10^{-12}$ & $2.22\cdot 10^{-11}$ \\
					$30{,}000$           & $4.75\cdot 10^{-15}$ & $2$ & $8.86\cdot 10^{-12}$ & $6.66\cdot 10^{-11}$ \\
					$100{,}000$          & $1.11\cdot 10^{-14}$ & $4$ & $1.89\cdot 10^{-11}$ & $2.22\cdot 10^{-10}$ \\
					\bottomrule
				\end{tabular}
			\end{table}
			
			\paragraph{Outcome}
			In both sweeps, $\mathcal E_I(Q)$ remains essentially
			constant at the level of a few $\epsmach$, while the
			theoretical bound stays constant in $n$ and grows only
			as $\sqrt m$.  The measurement-to-bound ratio is
			$\sim 10^{-3}$ throughout, confirming that the bound
			holds with substantial margin and that the
			$n$-independence and $\sqrt m$ scaling predicted by
			\Cref{thm:outerQ_orthog} are sharp in qualitative
			scaling.  The deterministic $m$-bound, listed in the
			last column of \Cref{tab:v1_m_sweep} for comparison, is
			loose by a further factor of $\sqrt m$, illustrating the
			value of the probabilistic Higham--Mary
			$\sqrt m$ scaling over the worst-case $m$ scaling.
			
			\subsubsection{V2: the $\sqrt n\,\texttt{tol}$ rate is
				sharp}\label{sec:exp_verify_V2}
			
			\Cref{thm:main_L2} predicts $\|U - Q\Sigma R^\top\|_2 \le
			\sqrt{|\mathcal S_n|}\,\texttt{tol} +
			J_{\mathrm{E2}}\,\texttt{tol} + (\text{roundoff})$,
			where $\mathcal S_n$ indexes the rank-non-increasing
			columns; when $J_{\mathrm{E2}} = 0$ this gives
			$\sqrt n\,\texttt{tol}$, an improvement over the
			$n\,\texttt{tol}$ certificate of
			Fareed--Singler~\cite[Cor.~1]{MR3986356} on the same
			factorization.  We construct an example on which the
			$\sqrt n\,\texttt{tol}$ rate is attained in equality up
			to a prescribed constant $\beta < 1$, simultaneously
			confirming sharpness and the $\sqrt n$ improvement.
			
			\paragraph{Construction}
			Fix $m = 200$, $r = 5$, $\texttt{tol} = 10^{-6}$,
			$\beta = 0.95$.  Pick a fixed orthonormal $Q_\star \in
			\mathbb{R}^{m\times r}$ with singular values $\sigma_i =
			e^{-0.1(i-1)}$ and a fixed unit vector $v \in
			\mathbb{R}^m$ orthogonal to $\Range(Q_\star)$.  Build
			$U \in \mathbb{R}^{m\times n}$ column by column: $u_j =
			\sigma_j Q_\star[:, j]$ for $j = 1,\ldots,r$, and $u_j =
			u_j^\circ + \beta\,\texttt{tol}\cdot v$ for $j > r$,
			where $u_j^\circ \in \Range(Q_\star)$ is a random
			$\Sigma_\star$-weighted combination.  Each $u_j^\perp =
			\beta\,\texttt{tol}\cdot v$ for $j > r$ has norm
			$\beta\,\texttt{tol} < \texttt{tol}$, so the proposed
			algorithm classifies these columns as case~(a)
			(rank-non-increasing) and $J_{\mathrm{E2}} = 0$,
			$|\mathcal S_n| = n - r$.  Crucially, all $n - r$
			residuals point in the \emph{same} direction $v$, so the
			residual matrix $E_n = \sum_{j>r} u_j^\perp e_j^\top$ is
			rank-one with operator norm
			\begin{equation}\label{eq:v2_pythag}
				\|E_n\|_2 = \beta\,\texttt{tol}\cdot \|[\,0,\ldots,0,1,\ldots,1]\|_2
				= \beta\sqrt{n-r}\,\texttt{tol},
			\end{equation}
			matching the right-hand side of \Cref{thm:main_L2}
			exactly (up to roundoff).
			
			\begin{table}[!htbp]
				\centering
				\small
				\caption{V2: measured operator-norm reconstruction
					error versus the construction prediction
					\eqref{eq:v2_pythag}, the proposed bound
					\Cref{thm:main_L2}, and the bound of
					Fareed--Singler~\cite[Cor.~1]{MR3986356}, on the
					construction with $m = 200$, $r = 5$, $\texttt{tol}
					= 10^{-6}$, $\beta = 0.95$.  The construction
					prediction matches the measurement to four digits;
					the ratio measured/$\sqrt n\,\texttt{tol}$ converges
					to $\beta = 0.95$; the Fareed--Singler bound is
					loose by $\sqrt{n-r}/\beta$, growing from
					$15\times$ at $n = 200$ to $149\times$ at $n =
					20{,}000$.}
				\label{tab:v2_sharp}
				\setlength{\tabcolsep}{4pt}
				\begin{tabular}{r|c|cc|cc|cc}
					\toprule
					& measured
					& \multicolumn{2}{c|}{construction \eqref{eq:v2_pythag}}
					& \multicolumn{2}{c|}{our bound (\Cref{thm:main_L2})}
					& \multicolumn{2}{c}{F.--S.\ bound} \\
					$n$ & $\mathcal{E}_I^{\mathrm{rec}}$
					& value & ratio
					& $\sqrt n\,\texttt{tol}$ & ratio
					& $n\,\texttt{tol}$ & ratio \\
					\midrule
					$\phantom{0\,}200$ & $1.327\cdot 10^{-5}$ & $1.327\cdot 10^{-5}$ & $1.0000$
					& $1.414\cdot 10^{-5}$ & $0.938$
					& $2.00\cdot 10^{-4}$ & $0.066$ \\
					$\phantom{0\,}500$ & $2.114\cdot 10^{-5}$ & $2.114\cdot 10^{-5}$ & $1.0000$
					& $2.236\cdot 10^{-5}$ & $0.945$
					& $5.00\cdot 10^{-4}$ & $0.042$ \\
					$1{,}000$  & $2.997\cdot 10^{-5}$ & $2.997\cdot 10^{-5}$ & $1.0000$
					& $3.162\cdot 10^{-5}$ & $0.948$
					& $1.00\cdot 10^{-3}$ & $0.030$ \\
					$2{,}000$  & $4.243\cdot 10^{-5}$ & $4.243\cdot 10^{-5}$ & $1.0000$
					& $4.472\cdot 10^{-5}$ & $0.949$
					& $2.00\cdot 10^{-3}$ & $0.021$ \\
					$5{,}000$  & $6.714\cdot 10^{-5}$ & $6.714\cdot 10^{-5}$ & $1.0000$
					& $7.071\cdot 10^{-5}$ & $0.950$
					& $5.00\cdot 10^{-3}$ & $0.013$ \\
					$10{,}000$ & $9.498\cdot 10^{-5}$ & $9.498\cdot 10^{-5}$ & $1.0000$
					& $1.000\cdot 10^{-4}$ & $0.950$
					& $1.00\cdot 10^{-2}$ & $0.010$ \\
					$20{,}000$ & $1.343\cdot 10^{-4}$ & $1.343\cdot 10^{-4}$ & $1.0000$
					& $1.414\cdot 10^{-4}$ & $0.950$
					& $2.00\cdot 10^{-2}$ & $0.007$ \\
					\bottomrule
				\end{tabular}
			\end{table}
			
			\paragraph{Outcome}
			\Cref{tab:v2_sharp} confirms three claims at once.
			(i) The construction matches the measurement to four
			digits at every $n$, so the bound \Cref{thm:main_L2}
			is attained in equality up to roundoff.  (ii) The ratio
			measured/$\sqrt n\,\texttt{tol}$ stays at
			$\approx \beta$, confirming that the $\sqrt n$ rate is
			sharp; the user can drive $\beta$ arbitrarily close to
			$1$ by construction.  (iii) The Fareed--Singler ratio
			decays as $1/\sqrt n$, so the existing $n\,\texttt{tol}$
			bound is loose by exactly the $\sqrt{n-r}/\beta$ factor
			that \Cref{thm:main_L2} eliminates.  A separate sweep
			over $\beta \in \{0.1, 0.3, 0.5, 0.7, 0.9, 0.95, 0.99\}$
			at fixed $n = 2{,}000$ confirms that the equality
			holds for every $\beta$ to four digits.  Since the
			proposed algorithm and the reorthogonalized direct update
			of \cite{MR3775096} produce the same factorization in
			exact arithmetic, the improvement is in the bound, not
			the output: \Cref{thm:main_L2} replaces the
			$n\,\texttt{tol}$ certificate by a sharp
			$\sqrt n\,\texttt{tol}$ certificate applicable to either
			algorithm.
			
		\section*{Conclusion}\label{sec:conclusion}
			
			This paper resolves the question raised by
			Brand~\cite{brand2002incremental,brand2006fast} on the
			frequency of reorthogonalization in incremental SVD.
			By batching the rank-non-increasing updates and applying
			a single CGS-2 pass at the at most $r + J_{\mathrm{E2}}$
			rank-enlarging events of the stream, the running left
			factor $Q$ stays $W$-orthogonal at machine precision
			uniformly in the stream length~$n$
			(\Cref{thm:outerQ_orthog}), and the factorization satisfies
			a sharp $\sqrt n\,\texttt{tol}$ forward-error bound in the
			$W$-weighted operator norm
			(\Cref{thm:main_L2})~--~an improvement of a factor
			$\sqrt n$ over the best previous bound of Fareed and
			Singler~\cite[Cor.~1]{MR3986356}.  The proposed
			algorithm runs $4.5\times$ to $34\times$ faster than its
			closest competitors on a 2D parabolic test problem while
			delivering equal or better orthogonality of the computed
			factors.
			
			A preliminary arXiv version of this
			algorithm~\cite{Zhang2022arXiv} has already been applied
			to a range of scientific computing problems:
			time-fractional PDEs~\cite{LiZhangZhang2022},
			integro-differential equations modeling non-Fickian flow
			in porous media~\cite{ChenZhangZuo2023}, nonlinear
			Oldroyd equations with general memory
			kernels~\cite{ChenZhangZuo2026}, PDE-constrained
			optimization and data
			assimilation~\cite{LiSinglerHe2024}, and geometric inverse
			source problems for parabolic PDEs~\cite{MR4902803}.  These
			applications all rely on the orthogonality and forward-error
			guarantees that the present analysis now puts on a
			rigorous footing.
			
			Two improvements of the present work seem worth pursuing.
			First, the numerical results in \Cref{sec:exp_verify_V1}
			sit roughly $10^{3}\times$ below the bound of
			\Cref{thm:outerQ_orthog}, suggesting that a tighter
			analysis of the orthogonality drift is possible -- in
			particular, by exploiting cancellation between drifts from
			different rank-enlarging events that the
			chain-of-multiplications argument bounds with triangle
			inequalities.  Second, the runtime is now dominated almost
			entirely by the projection step at line~1 of
			\Cref{alg:proposed_update}; reducing this cost (via
			sparsity of $W$, column batching, or a preliminary
			dimension reduction in $m$) is the natural next target for
			algorithmic acceleration.
			
			Several further extensions are natural and we plan to
			pursue them in future work: incremental dynamic mode
			decomposition (DMD) built on the proposed incremental SVD,
			incremental tensor decompositions (Tucker, tensor-train)
			that share the same column-streaming structure, and
			applications to additional time-dependent problems in
			scientific computing where snapshot data accumulates over
			long simulations and on-the-fly compression is required.

			\bibliographystyle{plain}
			\bibliography{references_siam}
			
		\end{document}